\documentclass[a4paper]{siamart220329} 

\usepackage{amstext,amsmath,amssymb}
\usepackage{epsfig,rotating,multirow}
\usepackage{graphicx,bm,color}
\usepackage{subcaption}
\usepackage{epstopdf}
\usepackage[noend]{algorithmic}
\usepackage{algorithm}
\usepackage{amsopn}
\floatstyle{plaintop}
\restylefloat{algorithm}
\usepackage{caption}
\newcommand{\bs}[1]{\boldsymbol{#1}}

\definecolor{gray}{gray}{0.6}

\def\bse{\bs{e}}
\def\bsf{\bs{f}}

\def\bsr{\bs{r}}
\def\bss{\bs{s}}
\def\bst{\bs{t}}
\def\bsu{\bs{u}}
\def\bsv{\bs{v}}
\def\bsw{\bs{w}}
\def\bsx{\bs{x}}
\def\bsy{\bs{y}}

\headers{Parallel-in-time MK}{Y. Erlangga}
\title{Parallel-in-time Multilevel Krylov Methods: A Prototype
}

\author{Yogi A. Erlangga\thanks{Zayed University, Mathematics and Statistics Department, Abu Dhabi Campus, P.O. Box 144534, Abu Dhabi, United Arab Emirates (\email{yogi.erlangga@zu.ac.ae, yogiae@gmail.com})}.}

\begin{document}

\maketitle

\begin{abstract}
This paper presents a parallel-in-time multilevel iterative method for solving differential algebraic equation, arising from a discretization of linear time-dependent partial differential equation. The core of the method is the multilevel Krylov method, introduced by Erlangga and Nabben~{\it [SIAM J. Sci. Comput., 30(2008), pp. 1572--1595]}. In the method, special time restriction and interpolation operators are proposed to coarsen the time grid and to map functions between fine and coarse time grids. The resulting Galerkin coarse-grid system can be interpreted as time integration of an equivalent differential algebraic equation associated with a larger time step and a modified $\theta$-scheme. A perturbed coarse time-grid matrix is used on the coarsest level to decouple the coarsest-level system, allowing full parallelization of the method. Within this framework, spatial coarsening can be included in a natural way, reducing further the size of the coarsest grid problem to solve. Numerical results are presented for the 1- and 2-dimensional heat equation using {\it simulated} parallel implementation, suggesting the potential computational speed-up of up to 9 relative to the single-processor implementation and the speed-up of about 3 compared to the sequential $\theta$-scheme.
\end{abstract}



\section{Introduction} \label{sec:intro}

In the early 2000's, Lions et al.~\cite{Lions01MT, Maday05T} introduced a time integration method, called {\it parareal}, for numerically solving the linear initial-boundary value problem:
\begin{eqnarray}
  \frac{\partial \widehat{u}}{\partial t} = f\left(\widehat{u},\frac{\partial \widehat{u}}{\partial x},\frac{\partial^2 \widehat{u}}{\partial x^2}\right),  \text{ in } (0,T] \times \Omega, \label{eq:pde_in_time}
\end{eqnarray}
with $x \in \Omega \subset \mathbb{R}^d$, the initial condition $\widehat{u}(0,x) = \widehat{u}_0$, and some boundary conditions. In constrast to the classical approaches that require sequential, time marching solves  starting from the initial condition/solution, parareal computes the solution over the time interval $\mathcal{T} = (0,T]$ in parallel, in an iterative fashion. 

To describe parareal, consider a spatially discretized version of~\eqref{eq:pde_in_time}:
\begin{eqnarray}
    \frac{\partial u(t)}{\partial t} = A_n u(t) + g(t),  \label{eq:dae}
\end{eqnarray}
where $A_n \in \mathbb{R}^{n \times n}$, with $n$ the number of spatial grid points (or nodal points), $u(t) \in \mathbb{R}^n$ is the solution vector at an instant time $t \in \mathcal{T}$,  and $g$ is the boundary condition vector. Applying the $\theta$-scheme for time integration results in the recurrence
$$
  u_k - u_{k-1} = (1-\theta) \Delta t A_n u_{k-1} + \theta \Delta t A_n u_k + \Delta t ((1-\theta) g_{k-1} + \theta g_k), 
$$
for  $k = 1,2,\dots,n_t$ and $\theta \in [0,1]$. Here, $u_k := u(t_k) = u(k \Delta t)$, where $\Delta t = T/n_t$ is the {\it fine} time step, used to compute the highly accurate solution in $\mathcal{T}$, with $n_t$ be the number of {\it fine} time steps. Note that $\theta = 0, \frac{1}{2}, 1$ corresponds to the explicit, Crank-Nicolson, and implicit scheme, respectively. Rearranging the above equation results in the recurrence
\begin{eqnarray}
   \psi_n u_k = \phi_n u_{k-1} + \Delta t \bar{g}_k, \label{eq:thetascheme}
\end{eqnarray}
where $\phi_n = I_n + (1-\theta) \Delta t A_n$, $\psi_n = I_n - \theta \Delta t A_n$, with $I_n$ be the identity matrix of size $n$, and $\bar{g}_k = (1-\theta) g_{k-1} + \theta g_k$. Starting from the initial solution $u_0 = u(0)$, the solution at any time $t_k = k \Delta t$ is computed sequentially using the recipe
\begin{align}
    u_k = \psi_n^{-1} (\phi_n u_{k-1} + \Delta t \bar{g}_k). \label{eq:timemarch}
\end{align}

Let $\mathcal{T}$ be partitioned into $n_T$ subintervals $\mathcal{T}_k = [T_{k-1},T_k]$, $k = 1,\dots,n_T$, with $|\mathcal{T}_k| = \Delta T = \nu \Delta t$, $1 < \nu \in \mathbb{N}$, $T_0 = 0$ and $T_{n_T} = T$. The {\it exact}\footnote{Here, the term``exact'' is referred to as the solution of the fully discretized counterpart of the initial-value problem, and not that of~\eqref{eq:pde_in_time}.} solution at the time $T_k = t_{\nu(k-1)} + \nu \Delta t$ is given by
\begin{align}
   U_k := U(T_k) = u_{\nu k} =   \Upsilon_{\nu} U_{k-1} + \tilde{g}_{\nu}, \quad k = 1,2, \dots,n_T, \label{eq:exactcoarse}
\end{align}
where $\Upsilon_{\nu} = (\psi^{-1}_n \phi_n)^{\nu}$, $\tilde{g}_{\nu} = \Delta t \psi^{-1}_n \displaystyle \sum_{m = 0}^{\nu-1} (\phi_n \psi^{-1}_n)^{m-\nu+1} g_{m + \nu(k-1)}$. This approach performs less time stepping, but requires access to the operator $\Upsilon_{\nu}$, which is not easily available. This computation can not also be performed in parallel because of the lack of initial condition at $T_{k-1}$ in the subinterval $\mathcal{T}_k$, except for $k = 1$.

To allow for parallelization, parareal replaces~\eqref{eq:exactcoarse} by the {\it coarse} time-step analogue to \eqref{eq:thetascheme}:
\begin{align}
     U_k = \Psi^{-1}_n (\Phi_n U_{k-1} + \Delta T \bar{G}_k), \label{eq:timemarchdT}
\end{align}
where  $\Phi_n = I_n + (1-\theta) \Delta T A_n$, $\Psi_n = I_n - \theta \Delta T A_n$, and $\bar{G}_k = (1-\theta) G_{k-1} + \theta G_k$, and uses it as an estimator for the initial condition in each subinterval $\mathcal{T}_k$. With an initial-value problem now defined in each subinterval $\mathcal{T}_k$, time integration using the fine time step $\Delta t$ based on~\eqref{eq:thetascheme} can be computed in parallel. This sequence of computation on the coarse and fine time grid however only produces the exact solution in $\mathcal{T}_1$, where the exact initial condition is used. A more accurate solution in the entire computational domain $\mathcal{T}$ is computed  by recalculating the coarse time-grid solution via a sequential correction procedure, which expectedly produces a more accurate initial condition for each subinterval $\mathcal{T}_k$. The new fine time solution in each time subinterval $\mathcal{T}_i$ is then updated using this improved initial condition, in parallel. This recipe, summarized in Algorithm~\ref{algo:parareal}, is repeated until a convergence on the fine time grid is achieved.

\begin{algorithm}[h!]
\caption{{\it Parareal}} 
\label{algo:parareal}
\begin{algorithmic}[1]
\STATE $U^0_0 = u_0$;
\FOR {$k = 1,\dots,n_T$}
        \STATE    $U^0_k = \Psi^{-1}_n (\Phi_n U^0_{k-1} + \Delta T \bar{G}_k)$;
\ENDFOR
\FOR {$\ell = 1,2,\dots$ until convergence }
    \FOR {$k = 1,\dots,n_T$}
          \FOR {$j = 1,\dots,\nu$}
                \STATE $u^{\ell}_{\nu (k-1) + j} = \psi^{-1}_n (\phi_n u^{\ell-1}_{\nu (k-1) + j -1} + \Delta t \bar{g}_{\nu(k-1)+j-1})$;
           \ENDFOR
   \ENDFOR
   \STATE $U^\ell_0 = U_{\ell-1}^0 = u_0$;
   \FOR {$k = 1,\dots,n_T$}
       \STATE $U^\ell_k = u_{\nu k}^0 +  \Psi^{-1} \Phi (U^{\ell}_{k-1} - U^{\ell-1}_{k-1})$;
   \ENDFOR
\ENDFOR
\end{algorithmic}
\end{algorithm}

Since the introduction of parareal, interests in time parallelization have arisen in the last two decades. Parareal, however, may not historically be regarded as the first work on time parallelization. The earliest work on this can be traced back to the work by Nievergelt~\cite{Nievergelt64} in the 1960s, an era where parallel computers did not even exist. While parareal and Nievergelt's method may look like two distant methods, they can in fact be classified under the same family: shooting methods; see~\cite{Gander15}. There are also methods which were developed in the context time-space domain decomposition and multigrid methods~\cite{Gander96, Gander03HN, Hoang13JJKR, Hackbusch84,Horton95V,Neumuller14}. We refer the reader to~\cite{Gander15} for a literature account on time parallelization methods prior to parareal.

There are some interpretations of the initial version of parareal, which algebraically links parareal with well-known methods for statiationary problems, and which has led to several alternative parallel-in-time approaches for solving initial-value problems~\cite{Gander07V, Gander15}. One interpretation relates parareal with a two-(time)grid method, with a nonstandard interpolation and restriction operator. Such an interpretation has allowed the construction of a parallel-in-time multigrid method, called MGRIT (multigrid reduction in time), where time parallelization is done in a multilevel framework~\cite{Falgout14FKMS}. At present, MGRIT is one of the most heavily tested and popular methods for time parallelization. Further development of MGRIT includes the inclusion of spatial coarsening in the method and extension to nonlinear evolution problems~\cite{Dobrev17KPS, Falgout17MOS, Falgout19LW, Howse19SFMS}.

In this paper, we propose another approach to time parallelization based on the multilevel Krylov (MK) framework, introduced in~\cite{Erlangga08Na}. We note here that, as a multilevel method, MK can not be considered as another multigrid method, due to, e.g., the nonexistence of smoothing process required for an effective multigrid method. In fact, MK as outlined in~\cite{Erlangga08Na} uses GMRES recursively to solve the coarse-grid system instead of multigrid smoothers such as Jacobi or Gauss-Seidel. 
Nevertheless, MK shares components that also are present in multigrid: fine-coarse grid maps and coarse-grid system. This allows the use of standard multigrid components in MK~\cite{Erlangga08Nc}. We shall demonstrate in this paper through numerical tests that a fast convergence of MK can still be achieved by using a simple technique for the fine-coarse time-grid mapping. Even though GMRES used on each coarse level increases the computational complexity of the method, we shall also show that for the model problem considered in the paper (the heat equation), the parallel-in-time MK method can be made a faster method than the serial $\theta$-scheme, when space coarsening is also performed. Due to the author's limited computational resources, the numerical tests were, however, only performed in a simulated parallel setting, where a parallel code is run sequentially on a single processor. 

The remainder of this paper is organized as follows. In Section~\ref{sec:mg-mkrylov}, to motivate the construction of our parallel-in-time multilevel Krylov method, we discuss the relation between parareal and two-grid methods. Section~\ref{sec:2level} presents the construction of the two-level Krylov method for time parallelization and discusses numerical tests and performance results using simulated parallel implementation. Extension to the multilevel Krylov method and numerical results are detailed in Section~\ref{sec:multilevel} together with numerical results. In Section~\ref{sec:conclusion}, we draw some concluding remarks and directions for future's work.

\section{Parareal as a two-grid method} \label{sec:mg-mkrylov}

To facilitate a link between parareal and a two-grid method, we first rewrite the recurrence~\eqref{eq:thetascheme} as a large system of equations
\begin{equation}
   \mathcal{A}_h \bsu_h = \bsf_h,  \label{eq:finesys}
\end{equation}
where
\begin{equation}
\mathcal{A}_h = 
\begin{bmatrix}
  I_n   &               &           &             \\
  -\phi_n & \psi_n  &           &      \\
          &   \ddots & \ddots &  \\
          &   &  -\phi_n  & \psi_n
\end{bmatrix},
\bsu_h = 
\begin{bmatrix}
  u_0  \\
  u_1  \\
  \vdots  \\
  u_{n_t}
\end{bmatrix}, \text{ and } \bsf_h =
\begin{bmatrix}
  u_0  \\
  \Delta t \bar{g}_1  \\
  \vdots  \\
  \Delta t \bar{g}_{n_t}
\end{bmatrix}.
\end{equation}
Here, $\mathcal{A}_h$ is an $(n_t+1) \times (n_t+1)$ block matrix. The serial solve of~\eqref{eq:thetascheme} is equivalent to solving \eqref{eq:finesys} using block forward substitution. We shall refer this system to as the fine time-grid system throughout.

Let the $(n_T +1)\times (n_T+1)$ block matrix
\begin{eqnarray}
   \mathcal{A}_H = \begin{bmatrix}
  I_n   &         &          &          \\
  -\Phi_n & \Psi_n   &          &          \\
          &   \ddots  &  \ddots &  \\
          &         &     -\Phi_n  & \Psi_n
\end{bmatrix} \label{eq:coarsemat}
\end{eqnarray}
be the coarse time-grid matrix analogue to $\mathcal{A}_h$. Furthermore, let $\bsu_H = [U_k]_{k=0}^{n_T}$, with $U_k = u_{\nu k}$, and $\bsf_H = [u_0^{\top},\Delta T \bar{G}_1^{\top}, \dots, \Delta T \bar{G}_{n_T}^{\top}]^{\top}$, with the superscript ``$^{\top}$'' denoting the transposition of a matrix of vector.

To reformulate parareal as a two-grid method, we define the following specific choice of components:
\begin{enumerate}
\item the prolongation matrix $\mathcal{P}$, defined as a $(n_t+1) \times (n_T+1)$ block matrix where the $(i,k)$-block is zero, except when $i = 1+\nu k$, $k = 0,\dots,n_T$, where the block $\mathcal{P}_{i,k} = I_n$;
\item the restriction matrix $\mathcal{R}$, defined as a $(n_T+1) \times (n_t+1) $ block matrix where the $(k,i)$-block is zero, except when $i = 1+\nu k$, $k = 0,\dots,n_T$, where the block $\mathcal{R}_{i,k} = \Phi_n^{-1} \psi_n$;
\item the coarse time-grid matrix $\mathcal{A}_H$, given by~\eqref{eq:coarsemat}. Note that this coarse time-grid matrix is not of Galerkin type, i.e.,  $\mathcal{A}_H \neq \mathcal{R} \mathcal{A}_h \mathcal{P}$;
\item smoother, which is represented algebraically by the $n_T \times n_T$ modified block Jacobi matrix
\begin{equation}
   \mathcal{M}_{J} = \begin{bmatrix}
                   \bar{M}_{J} &             &           &            \\ 
                                        & M_{J} &           &            \\
                                        &             & \ddots &            \\
                                        &             &           & M_{J} 
          \end{bmatrix},
\end{equation}
where the diagonal blocks
\begin{equation}
    \bar{M}_{J} = \begin{bmatrix}
           I_n  &          &            & \\
          -\phi_n &  \psi_n   &           & \\ 
                  & \ddots & \ddots & \\
                  &          &   -\phi_n   & \psi_n 
   \end{bmatrix} \text{ and } M_{J} = \begin{bmatrix}
                  \psi_n   &           &            & \\ 
                 -\phi_n  & \ddots &            &  \\
                          & \ddots & \ddots  & \\
                          &           & -\phi_n    & \psi_n
   \end{bmatrix}
\end{equation}
are an $(\nu+1) \times (\nu+1)$ and  $\nu \times \nu$ block matrix, respectively. Furthermore, let $\mathcal{I}_{J}$ be the $n_t \times n_t$ modified block identity matrix with zero matrix on the block that corresponds to the $u_{\nu k}$ unknown, for $k = 0,1,\dots, n_T$.
\end{enumerate}

It can be shown that the following two-grid procedure produces parareal iterands $u_k^{\ell}$ and $U_k^{\ell}$:
\newline
\begin{enumerate}
    \item Set the initial approximate solution $\bsu^0_h$ using the following steps:
            \begin{enumerate}
                  \item[1.1.] Solve the coarse time-grid problem:
                                   $$
                                         \mathcal{A}_H \bsu_H^0 = \bsf_H.
                                   $$
                 \item[1.2.] Prolongate the coarse time-grid solution to the fine time-grid:
                                   $$
                                        \bsu_h ^0 = \mathcal{P} \bsu^0_H.
                                   $$
            \end{enumerate}
     \item For $\ell = 1,2, \dots$ until convergence:
            \begin{enumerate}
                 \item[2.1.] Pre-smoothing:
                                 $$
                                   \widetilde{\bs{u}}_h^{\ell-1} = \bs{u}^{\ell-1}_h + \mathcal{I}_{J} \mathcal{M}^{-1}_{J}(\bs{f}_h - \mathcal{A}_h \bsu_h^{\ell-1}).
                                 $$
                 \item[2.2.] Restriction:
                                \begin{align}
                                     \widetilde{\bs{r}}_H^{\ell-1} &= \mathcal{R} \left( \bs{f}_h - \mathcal{A}_h \widetilde{\bs{u}}_h^{\ell-1} \right). \notag
                                \end{align} 
               \item[2.3.] Coarse time-grid solve:
                          $$\bs{d}_H^{\ell-1} = \mathcal{A}_H^{-1} \widetilde{\bs{r}}_H^{\ell-1}.$$
               \item[2.4.] Prolongation and correction:
                          $$
                               \bs{u}_h^{\ell} = \widetilde{\bs{u}}_h^{\ell-1} + \mathcal{P} \bs{d}_H^{\ell-1}.
                          $$
            \end{enumerate}
\end{enumerate}

The above two-grid method produces a sequence of error vectors $\bse_h^{\ell}$ that satisfies the relation
\begin{equation}
    \bse_h^{\ell} = \mathcal{K}_{MG}  \bse_h^{\ell-1},  \label{eq:errorMG}
\end{equation}
where $\mathcal{K}_{MG} := (\mathcal{I}_h - \mathcal{P} \mathcal{A}_H^{-1} \mathcal{R} \mathcal{A}_h) (\mathcal{I}_h - \mathcal{I}_J \mathcal{M}_J^{-1} \mathcal{A}_h)$. 

In $\mathcal{K}_{MG}$, the term $\mathcal{I}_h - \mathcal{P} \mathcal{A}_H^{-1} \mathcal{R} \mathcal{A}_h  =: Q_{D,h}$, which represents the coarse-grid correction,  is similar to the deflation operator~\cite{Nico87,Frank01V,Erlangga08Nb}, used to accelerate a Krylov iterative method. However, like in standard multigrid methods, deflation uses Galerkin coarse-grid matrix and some conditions for $\mathcal{R}$ and $\mathcal{P}$ such that the Galerkin coarse-grid matrix is nonsingular. In this case, applying $\mathcal{Q}_{D,h}$ on $\mathcal{A}_h$ deflates or shifts some eigenvalues of $\mathcal{A}_h$ to zero, resulting in a convergence acceleration of a Krylov method~\cite{Nabben04V}. The number of eigenvalues of $\mathcal{A}_h$ deflated to zeros typically depends on the rank condition of the matrix $\mathcal{R}$ and $\mathcal{P}$.

Due to a non-Galerkin coarse-grid approach, $Q_{D,h}$ in this case does not play the role one typically expects from deflation: that near-zero eigenvalues in the spectrum of $\mathcal{A}_h$ are deflated so that they are no longer present in the spectrum of $\mathcal{A}_h Q_{D,h}$, leading to improved convergence. 

In~\cite{Tang09NVE}, several other variants of shifting-based preconditioning operators derived from multigrid and domain decomposition methods are derived and analyzed.

\section{Two-level Krylov for time parallelization} \label{sec:2level}

\subsection{General approach} 
Consider again the partition of the time interval $\mathcal{T} = [0,T]$ into $n_T$ uniform subintervals, $\mathcal{T}_k = [T_{k-1},T_k]$, $k = 1,\dots,n_T$, with $T_0 = t_0 = 0$ and $T_{n_T} = t_{n_t} = T$, such that $\Delta T = T/n_T = \nu T/n_t = \nu \Delta t$, $1 < \nu \in \mathbb{N}$. Using this partitioning, the linear system~\eqref{eq:finesys} can be rewritten as
\begin{align}
  \begin{bmatrix}
      \mathcal{A}^h_{0,0} &                           &  & \\
     \mathcal{A}^h_{1,0 }& \mathcal{A}^h_{1,1} &  & \\
                                & \ddots                  & \ddots &  \\
                                &                            & \mathcal{A}^h_{n_T,n_T-1} & \mathcal{A}^h_{n_T,n_T} 
  \end{bmatrix}  \begin{bmatrix}
                             \bsu^h_0 \\
                             \bsu^h_1 \\
                             \vdots \\
                             \bsu^h_{n_T}
                        \end{bmatrix} =  \begin{bmatrix}
                             \bsf^h_0 \\
                             \bsf^h_1 \\
                             \vdots \\
                             \bsf^h_{n_T}
                        \end{bmatrix}, \label{eq:blocksys}
\end{align}
where $\mathcal{A}_{0,0}^h = I_n$, $\bsu_0^h = u_0$, and, for $k = 1,\dots, n_T$,  $\bsu^h_k = [u^h_{\nu(k-1) + i}]_{i = 1,\dots,\nu}$ and $\mathcal{A}^h_{k,k}$ is a $\nu \times \nu$-block matrix.

The two-level Krylov method~\cite{Erlangga08Na} for~\eqref{eq:blocksys} is a Krylov method applied to the preconditioned system
\begin{eqnarray} \label{eq:AQsystem}
    \mathcal{A}_h \mathcal{Q}_h \widetilde{\bsu}_h = \bsf_h,
\end{eqnarray}
where
\begin{eqnarray}
 \mathcal{Q}_h = \mathcal{I}_h - \mathcal{Z} \mathcal{A}_H^{-1} \mathcal{Y}^{\top} \mathcal{A}_h +\mu \mathcal{Z} \mathcal{A}_H^{-1} \mathcal{Y}^{\top}, \quad \mathcal{A}_H = \mathcal{Y}^{\top} \mathcal{A}_h \mathcal{Z},  \label{eq:Qh}
\end{eqnarray}
with $\mu \in \mathbb{R}$. Here, both $\mathcal{A}_H$ and $\mathcal{Q}_h$ are assumed to be nonsingular so that $\widetilde{\bsu}_h = \mathcal{Q}_h^{-1} \bsu^h$ is defined. 

The choice $\mu = 0$ in~\eqref{eq:Qh} corresponds to the deflation operator~\cite{Frank01V, Tang09NVE}. This choice however leads to a singular $\mathcal{Q}_h$ and is known to be not a good choice for a multilevel implementation, due to ``reappearance'' of near-zero eigenvalues in the spectrum of $\mathcal{A}_h \mathcal{Q}_h$ under an inaccurate inversion of the coarse time-grid matrix $\mathcal{A}_H$. Since near-zero eigenvalues generally are responsible for slow convergence, an effective deflation method requires an exact solve of the coarse time-grid system. Such a requirement limits the deflation approach not only to a two-level setting -- unless the coarse time-level system is solved very accurately using a multilevel technique -- but also the size of the coarse time-level system to solve.  

Alternatively, as used in~\cite{Erlangga08Na}, the parameter $\mu$ is set equal to the largest eigenvalues of $\mathcal{A}_h$; $\mu = \lambda_{\max}(\mathcal{A}_h)$. With this choice, the action of $\mathcal{Q}_h$ ``shifts'' some near-zero eigenvalues of $\mathcal{A}_h$ to the $\lambda_{\max}(\mathcal{A}_h)$. Under an inaccurate inversion of $\mathcal{A}_h$, some eigenvalues of $\mathcal{A}_h \mathcal{Q}_h$ are expected to appear slightly away from $\lambda_{\max}(\mathcal{A}_h)$, but remain far from zero. This property allows the use of an inaccurate inversion of $\mathcal{A}_H$ using, for instance, an iterative method, making the method practical. Furthermore, with the spectrum of $\mathcal{A}_h \mathcal{Q}_h$ distributed away from zero, we may expect an improved convergence of a Krylov method.

The matrices $\mathcal{Y}, \mathcal{Z} \in \mathbb{R}^{n_t n \times r}$ in~\eqref{eq:Qh} are called the {\it deflation} matrices, assumed to be of rank $r < n_t n$ and such that $\mathcal{A}_H = \mathcal{Y}^T \mathcal{A}_h \mathcal{Z}$ is nonsingular. From a theoretical view point, $\mathcal{Y}$ and $\mathcal{Z}$ are chosen such that the spectrum of $\mathcal{A}_h \mathcal{Q}_h$ has no near-zero eigenvalues. An example of such a choice is when the columns of $\mathcal{Y}$ and $\mathcal{Z}$ are left and, respectively, right eigenvectors associated with the $r$ near-zero eigenvalues of $\mathcal{A}_h$. In this case, $r$ near-zero eigenvalues are shifted to $\lambda_{\max}$, with the other $n_t n-r$ eigenvalues remaining untouched. Furthermore, a larger $r$ means more near-zero eigenvalues are shifted, resulting in a more clustered eigenvalues of $\mathcal{A}_h\mathcal{Q}_h$, and expectedly faster convergence. Computing a set of $r$ eigenvectors is however very costly, especially when the size of $\mathcal{A}_h$ is large. Since $\mathcal{Y}$ and $\mathcal{Z}$ algebraically resemble integrid operators in multigrid (cf. $\mathcal{K}_{MG}$ below Eq.~\eqref{eq:errorMG}), in practice, we shall use the easier-to-construct restriction and interpolation (or integrid) matrix known in multigrid.

A two-level Krylov implementation is presented in Algorithm~\ref{algo:MK}, where flexible GMRES (FGMRES)~\cite{Saad93} is used as the basis. 
\begin{algorithm}[h!]
\caption{{\it FGMRES-based Multilevel Krylov implementation}} 
\label{algo:MK}
\begin{algorithmic}[1]
\STATE Choose $\bsu_h^{0}$ and $\mu$.
\STATE Compute $\bsr_h^{0} = \bsf_h - \mathcal{A}_h \bsu_h^{0}$, $\beta = \| \bsr_h^0 \|_2$, and $\bsv_h^{1} = \bsr_h^{0}/\beta$.
\FOR {$j = 1,2,\dots, J$ }
   \STATE $\bsx_h^{j} := \mathcal{Q}_h \bsv_h^{j}$;
   \STATE $\bsw_h := \mathcal{A}_h \bsx_h^{j}$;
   \FOR {$l = 1,\dots,j$}
        \STATE $h_{l,j} := (\bsw_h, \bsv_h^{j})$;
        \STATE $\bsw_h :=  \bsw_h - h_{l,j} \bsv_h^j$;
   \ENDFOR
   \STATE $h_{j+1,j} := \|\bsw_h\|_2$;
   \STATE $\bsv_h^{j+1} := \bsw_h/h_{j+1,j}$;
\ENDFOR
\STATE Set $\mathcal{X}_J = [\bsx_h^{1} \dots \bsx_h^{J}]$ and $\hat{H}_J = [h_{i,j}] _{1\le i \le j+1; 1 \le j \le J}$;
\STATE Compute $\bsy_{J} = \arg \min_{\bsy} \|\beta \bse_1 - \hat{H}_J \bsy \|_2$ and $\bsu_h^J = \bsu_h^{0} + \mathcal{X}_J \bsy_{J}$.
\end{algorithmic}
\end{algorithm}
In Algorithm~\ref{algo:MK}, Step 4 corresponds to the (right-preconditioned) two-level Krylov step. In this step, the matrix $\mathcal{Q}_h$ is never explicitly formed. Instead, its action on the vector $\bsv_h^{j}$ is performed implicitly using the following substeps:
\begin{align}
   &\text{\small 4.1:}~\bss_h := \mathcal{A}_h \bsv_h^{j};& &\notag \\
   &\text{\small 4.2:}~\widehat{\bsv}_H := \mathcal{Y}^{\top}(\bss_h - \mu \bsv_h^{j});& &\notag \\
   &\text{\small 4.3:}~\widetilde{\bsv}_H := \mathcal{A}_H^{-1} \widehat{\bsv}_H;& & \notag \\
   &\text{\small 4.4:}~\bst_h := \mathcal{Z}\widetilde{\bsv}_H;& &\notag  \\
   &\text{\small 4.5:}~ \bsx_h^{j} := \bsv_h^{j} - \bst_h;& &\notag 
\end{align}


\subsection{Intergrid matrices}

To construct the matrix $\mathcal{Y}$ and $\mathcal{Z}$,  we shall repurpose the multigrid intergrid matrices. For instance, following MGRIT~\cite{Falgout14FKMS}, we can set $\mathcal{Y}$ as a matrix that corresponds to a direct injection process of a fine time-grid quantity to the coarse time-grid, i.e.,
\begin{align} \label{eq:mgritY}
   \mathcal{Y} = \begin{bmatrix}
                             0 \\
                          \mathcal{I}_c
                         \end{bmatrix} 
\end{align}
and $\mathcal{Z}$ the {\it ideal} interpolation:
\begin{align}  \label{eq:mgritZ}
   \mathcal{Z} = \begin{bmatrix}
                            - ( \mathcal{A}^h_{ff})^{-1} \mathcal{A}^h_{fc} \\
                              \mathcal{I}_{c}
                         \end{bmatrix}, 
\end{align}
where $\mathcal{A}^h_{ff}$ and $\mathcal{A}^h_{fc}$ are submatrices of $\mathcal{A}_h$ associated with partitioning  of unknowns into the fine time-grid and course time-grid unknowns:
$$
   \mathcal{A}_h = \begin{bmatrix}
                              \mathcal{A}_{ff}^h &  \mathcal{A}_{fc}^h \\
                              \mathcal{A}_{cf}^h &  \mathcal{A}_{cc}^h
                         \end{bmatrix}.
$$
Consequently, $\mathcal{A}_H = \mathcal{Y}^{\top} \mathcal{A}_h \mathcal{Z} =  \mathcal{A}_{cc}^h - \mathcal{A}_{cf}^h (\mathcal{A}_{ff}^h)^{-1}  \mathcal{A}_{fc}^h$. As in parareal, this expensive-to-invert $\mathcal{A}_H$ can be replaced by the coarse time-grid analogue of $\mathcal{A}_h$, associated with a sequential solve of time marching problem with the larger time step $\Delta T$, given by~\eqref{eq:coarsemat}.

For the  two-level Krylov method, we however consider the specific construction of the intergrid operators based on the $(n_T+1) \times (n_T+1)$-block diagonal matrix
\begin{eqnarray}
  \quad \mathcal{Y} = blockdiag(\mathcal{Y}_k), \quad \mathcal{Z} = blockdiag(\mathcal{Z}_k) \label{eq:blockYZ},
\end{eqnarray}
where the submatrices $\mathcal{Y}_k$ and $\mathcal{Z}_k$, $k = 0,\dots,n_T$ are full-rank. The size of $\mathcal{Y}_k$ and $\mathcal{Z}_k$ is assumed to be compatible under left- and right-multiplication with $\mathcal{A}^h_{k,k}$, so that the coarse time-grid matrix $\mathcal{A}_H$ is defined. In this case, 
\begin{align}
  \mathcal{A}_H = \begin{bmatrix}
                                  \mathcal{A}^H_{0,0} &                             &    &\\
                                  \mathcal{A}^H_{1,0} & \mathcal{A}^H_{1,1} &    &\\
                                                              & \ddots                   & \ddots & \\
                                                             &                              & \mathcal{A}^H_{n_T,n_T-1} & \mathcal{A}^H_{n_T,n_T}
                            \end{bmatrix}, \label{eq:coarsesys}
\end{align}
where the diagonal blocks $\mathcal{A}^H_{k,k} = \mathcal{Y}_{k}^{\top} \mathcal{A}^h_{k,k} \mathcal{Z}_k$, $ k = 0,\dots,n_T$, are nonsingular, and  $\mathcal{A}^H_{k,k-1} = \mathcal{Y}_{k}^{\top} \mathcal{A}^h_{k,k-1} \mathcal{Z}_{k-1}$, $k = 1,\dots,n_T$.


In the sequel, we discuss two constructions, which lead to time coarsening and time-space coarsening two-level Krylov methods. 

\subsubsection{Case 1: Time coarsening}

One possible setting for the intergrid matrices, which corresponds to the (pure) time  (or T-coarsening, throughout), is based on $\mathcal{Z}_0 = I_n$ and
\begin{align}  \label{eq:tcoarseZ}
   \mathcal{Z}_k = [I_n \, I_n \, \dots \, I_n]^{\top}, \quad k = 1,\dots,n_T,
\end{align}
where $\mathcal{Z}_k$ is a $\nu \times 1$ block matrix, and the setting $\mathcal{Y}_k = \mathcal{Z}_k$. Submatrices in the coarse-grid matrix $\mathcal{A}_H$~\eqref{eq:coarsesys} can be computed explicitly and are given by
\begin{align}
    \mathcal{A}^H_{0,0} &= I_n, \label{eq:timecoarse1}\\
    \mathcal{A}^H_{kk} &= \nu \psi_n - (\nu - 1)\phi_n  = I_n - (\nu - \theta)\Delta t A_n = I_n - (1-\tilde{\theta})\Delta T A_n =: \widetilde{\Psi}_n, \label{eq:timecoarse2}\\
   \mathcal{A}^H_{k,k-1} &= -\phi_n = -(I_n + \theta \Delta t A_n) = -(I_n + \tilde{\theta} \Delta T A_n) =: - \widetilde{\Phi}_n, \label{eq:timecoarse3}
\end{align}
for $k = 1,\dots,n_T$, where  $\tilde{\theta} = \theta/\nu \in [0,1/2]$. 

\begin{figure}[!h]
\begin{center}
\includegraphics[width=0.75\textwidth]{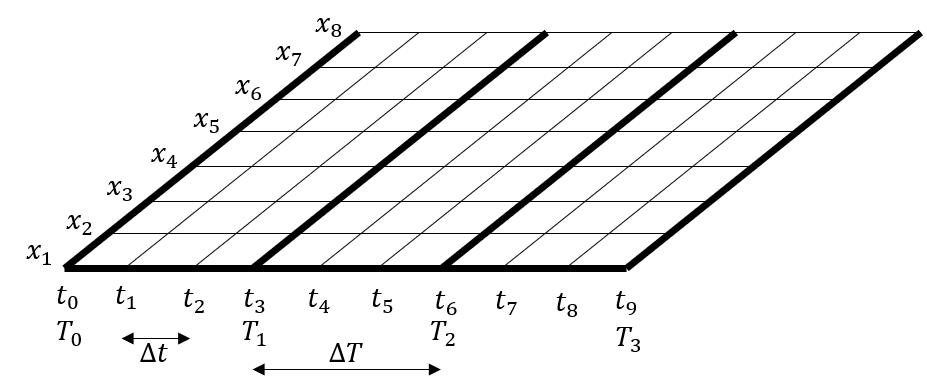}
\end{center}
\caption{Time coarsening, with $\Delta T = \nu \Delta t$, with $\nu = 3$.} \label{fig:time-coarsening}
\end{figure}

Note that the equations~\eqref{eq:timecoarse1}--\eqref{eq:timecoarse3} differ slightly from parareal's and MGRIT's coarse time-grid problem in the parameter $\theta$ (cf. Eq.~\eqref{eq:coarsemat}). With \eqref{eq:timecoarse1}--\eqref{eq:timecoarse3}, the coarse time-grid system can be interpreted as a modfied $\theta$-scheme with a larger time step $\Delta T = \nu \Delta t$ (see Figure~\ref{fig:time-coarsening}). A meaningful time coarsening requires $\nu \le 2$, resulting in $\tilde{\theta} < \theta$, which may lead to stability issues on the coarse grid for large $\Delta T$. The case $ \tilde{\theta}=\theta$ (as in parareal) only occurs when $\theta = 0$ (the explicit method).

\subsubsection{Case 2: Time-space coarsening}

With the computational work dominated by the inversion of the implicit term $\psi_n$, the total computational work of the sequential $\theta$-scheme is of order $n_t C_{inv}(n)$, with $C_{inv}$  the complexity associated with the action of $\psi_n^{-1}$ on a vector (for 2D cases, typically $C^h_{inv}(n) \sim n^2$). For the GMRES-based two-level Krylov method, the order of  computational work is 
$$
W_{2L} \sim N_{iter} \frac{n_t}{\nu} C^H_{inv},
$$
where $N_{iter}$ is the number of GMRES iterations needed to reach convergence and $C^H_{inv}$ is the complexity associated with the action of $\tilde{\Phi}_n^{-1}$ in the coarse time-grid system, which is the same as $\psi_n^{-1}$. Therefore, for a two-level method to be competitive, the ratio $N_{iter}/\nu$ has to be made less than 1, via an aggressive time coarsening with large $\nu$. This is however not easy to attain solely using an aggresive time coarsening. As we shall see in the numerical example, a large $\nu$ used in an aggressive time coarsening slows down the convergence of GMRES (so, $N_{iter}/\nu > 1$) that even makes the method impractical, due to a significant increase in the number of Arnoldi vectors to store.

A reduction in the computational complexity can be attained via a combination of time and space (TS) coarsening, resulting in a much smaller coarse time-space-grid system. For the construction of the time-space coarsening, we consider the intergrid matrices~\eqref{eq:blockYZ}, but now with
\begin{align}
   \mathcal{Z}_k = [Z \, Z \, \dots \, Z]^{\top}  \quad \text{and}  \quad \mathcal{Y}_k = [Y \, Y \, \dots \, Y]^{\top}, \label{eq:TS-coarsening}
\end{align}
where $\mathcal{Z}_k, \mathcal{Y}_k$ are a $\nu \times 1$-block diagonal matrix, with a full rank $Z \in \mathbb{R}^{n\times m}$, $m < n$. In this case, the submatrices in $\mathcal{A}_H$ are given by
\begin{align}
   \mathcal{A}^H_{0,0} &= Y^{\top} Z \in \mathbb{R}^{m \times m}, \\
    \mathcal{A}^H_{kk} &= Y^{\top}Z  - (1-\tilde{\theta})\Delta T Y^{\top} A Z \in \mathbb{R}^{m \times m}, \\
   \mathcal{A}^H_{k,k-1} &= -(Y^{\top}Z + \tilde{\theta} \Delta T Y^{\top}AZ) \in \mathbb{R}^{m \times m},
\end{align}
for $k = 1,\dots,n_T$. 

\begin{figure}[!h]
\begin{center}
\includegraphics[width=0.75\textwidth]{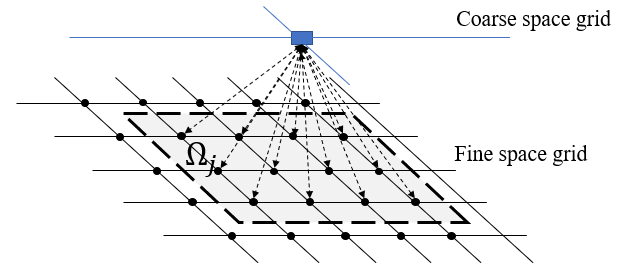}
\end{center}
\caption{Agglomeration of 12 fine spatial gridpoints in $\Omega_j \subset \Omega \subset \mathbb{R}^2$ into 1 coarse spatial gridpoint.} \label{fig:2d_space_coarsening}
\end{figure}

For $Y^{\top}$ and $Z$, we use  simple-to-construct agglomeration matrices, which are not a typical choice for multigrid but work very well for multilevel Krylov~\cite{Erlangga08Na}. We consider a partition of the spatial domain $\Omega \subset \mathbb{R}^d$ into $N$ nonoverlapping subdomains $\Omega_j$, such that $\Omega = \bigcup_{j=1}^N \Omega_j$. Let $\bsx_i \in \Omega$ be a spatial gridpoint, and let $\mathcal{I} = \{i | \bsx_i \in \Omega\}$ be the index set. Furthermore, let $\mathcal{I}_j = \{i \in \mathcal{I} | \bsx_i \in \Omega_j\}$. Figure~\ref{fig:2d_space_coarsening} illustrates this partitioning for $\Omega \subset \mathbb{R}^2$, with $\Omega_j \subset \Omega$ containing 12 gridpoints with indices in $\mathcal{I}_j$.

The matrix $Z = [z_{ij}]$ is defined as follows:
\begin{eqnarray}
z_{ij} = 
\begin{cases}
1,& i \in \mathcal{I}_j \\
0,& i \notin \mathcal{I}_j
\end{cases}.  \label{eq:constructZ}
\end{eqnarray}
The matrix $Y = [y_{ji}]$ is given by
\begin{eqnarray}
y_{ji} = 
\begin{cases}
1/|\mathcal{I}_j|,& i \in \mathcal{I}_j \\
0,& i \notin \mathcal{I}_j
\end{cases}. \label{eq:constructY}
\end{eqnarray}

\begin{figure}[!h]
\begin{center}
\includegraphics[width=0.85\textwidth]{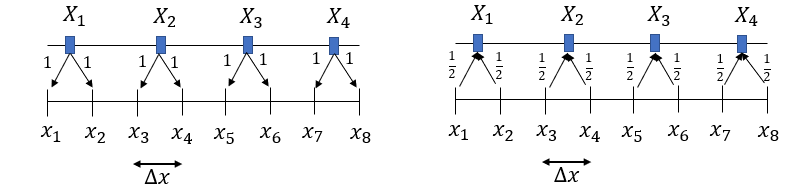}
\end{center}
\caption{1D agglomeration using standard spatial coarsening. Left: the construction of $Z$. Right: the construction of $Y$.} \label{fig:1d_space_coarsening}
\end{figure}

A specific example in 1D is shown in Figure~\ref{fig:1d_space_coarsening}, where standard spatial coarsening, i.e., where the coarse space grid is generated by agglomerating 2 fine spatial gridpoints, is used. One can show that $Y^{\top}Z = I_m$, with $m = 4$ for this particular example. Furthermore, $Y^{\top} AZ$ is associated with a direct discretization with the coarse grid size $\Delta X = \sqrt{2}\Delta x$.

\subsection{Two-level convergence} \label{sec:two-level}

To have insights on the convergence of the two-level Krylov method, we shall consider the heat equation
\begin{eqnarray}
     \frac{\partial \hat{u}}{\partial t} = \tau \sum_{i=1}^d \frac{\partial^2 \hat{u}}{\partial x_i^2}, \text{ in } \Omega \subset \mathbb{R}^d \quad (d \in \{1,2\}),
\end{eqnarray}
with the initial condition $u(0,\cdot) = u_0$ and the boundary conditions $u(\cdot,\bsx) = 0$, for $\bsx \in \partial \Omega$.
The spatial derivatives are discretized by the second-order finite difference scheme. The resultant differential algebraic equation is integrated using the $\theta$-scheme, with $\theta = 0.5$ (Crank-Nicolson). 

We shall first discuss spectra of the two-level matrix $\mathcal{A}_h\mathcal{Q}_h$ for the 1-dimensional case. By a direct computation it can be shown that
\begin{align}
   \mathcal{A}_h\mathcal{Q}_h = \begin{bmatrix}
                                                    \mathcal{A}^h_{0,0} \mathcal{Q}^h_0 &                &  & \\
                                                         \star   & \mathcal{A}^h_{1,1} \mathcal{Q}^h_1 &  & \\
                                                         \star   & \star         & \ddots  &  \\
                                                          \star   & \star         & \star  & \mathcal{A}^h_{n_T,n_T} \mathcal{Q}^h_{n_T} 
                                               \end{bmatrix},
\end{align}
where $\mathcal{Q}^h_k = \mathcal{I}^h_k - \mathcal{Z}_k (\mathcal{A}^H_{k,k})^{-1} \mathcal{Y}^{\top}_k \mathcal{A}^h_{k,k} + \mu \mathcal{Z}_k (\mathcal{A}^H_{k,k})^{-1} \mathcal{Y}^{\top}_k$. Consequently,
\begin{eqnarray}
      \sigma(\mathcal{A}_h \mathcal{Q}_h) = \bigcup_{k=0}^{n_T} \sigma(\mathcal{A}^h_{k,k} \mathcal{Q}^h_k),
\end{eqnarray}
where $\sigma(\cdot)$ denotes the spectrum of the matrix in the argument. Thus, the eigenvalues can be computed block-wise. Furthermore, for uniform coarsening (with fixed $\Delta T$), the block matrix $\mathcal{A}^h_{k,k} \mathcal{Q}^h_{k,k}$, $k = 1,\dots,n_T$, is constant. Thus, the complete spectrum can be obtained by computing eigenvalues of $\mathcal{A}^h_{k,k} \mathcal{Q}^h_{k,k}$, $k=0,1$. 

\subsubsection{1-dimensional case} \label{sec:1d-case}

\begin{figure}[!h]
\begin{center}
\begin{subfigure}{.32\textwidth}
\begin{center}
\includegraphics[width=1\textwidth]{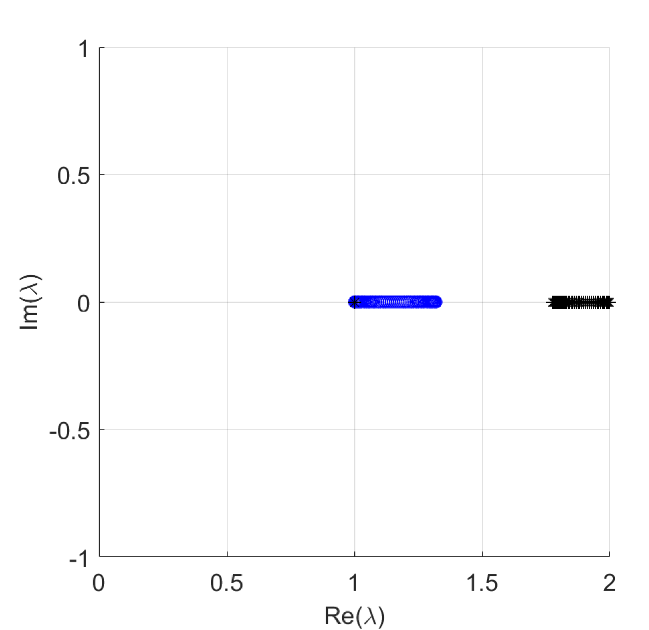}
\end{center}
\caption{$\nu = 2$}
\end{subfigure}
\begin{subfigure}{.32\textwidth}
\begin{center}
\includegraphics[width=1\textwidth]{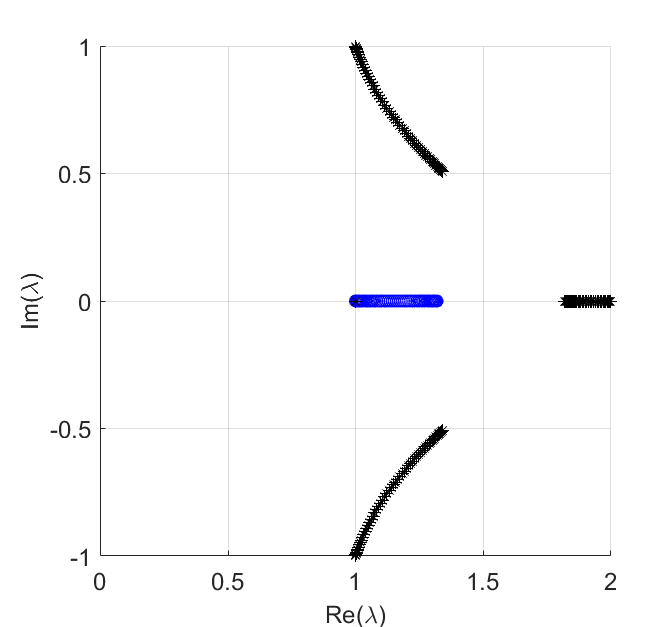}
\end{center}
\caption{$\nu = 4$}
\end{subfigure}
\begin{subfigure}{.32\textwidth}
\begin{center}
\includegraphics[width=1\textwidth]{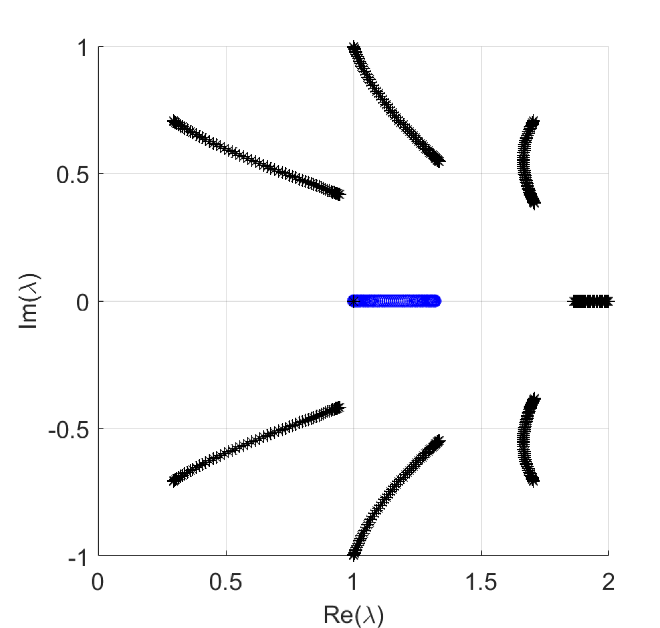}
\end{center}
\caption{$\nu = 8$}
\end{subfigure}
\end{center}
\caption{Spectrum of $\mathcal{A}_h \mathcal{Q}_h$ (in black), from the 1-dimensional heat equation. The two-level setting uses  time coarsening intergrid operators~\eqref{eq:tcoarseZ} and $\mathcal{Y}_k = \mathcal{Z}_k$. The spectrum of $\mathcal{A}_h$ is plotted in blue. For the discretization, $n_x = 63$, $n_t = 32$, the Courant number $C = 2k \Delta t/\Delta x^2 = 0.32$, and $\theta = 0.5$.} \label{fig:spectrum_1d_aggloA}
\end{figure}

Figure~\ref{fig:spectrum_1d_aggloA} shows computed eigenvalues (in black) of the two-level matrix $\mathcal{A}_h \mathcal{Q}_h$ using the T-coarsening for the time coarsening parameter $\nu \in \{2, 4, 8\}$ and the shift parameter $\mu = 1$. The spectrum of the original matrix $\mathcal{A}_h$, shown in blue in Figure~\ref{fig:spectrum_1d_aggloA}, lies along the real axis between 1 and 1.5. For $\nu = 2$, the spectrum of $\mathcal{A}_h \mathcal{Q}_h$ splits into two clusters: one at 1, which corresponds to the shift $\mu = 1$, and the other cluster spread between 1.5 and 2 in the real axis. For $\nu = 4$ and $8$, the spectra are generally less clustered with eigenvalues split in 4 and 8 clusters in the complex plane, respectively, with one cluster exactly at 1. With less clustered eigenvalues, we should expect slower convergence rate of the two-level method, when a more aggresive time coarsening approach is used.

\begin{figure}[!h]
\begin{center}
\begin{subfigure}{.32\textwidth}
\begin{center}
\includegraphics[width=1\textwidth]{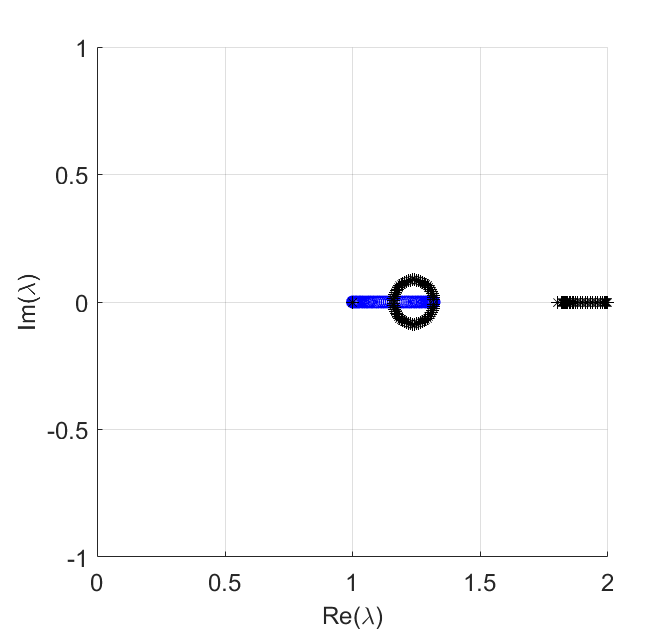}
\end{center}
\caption{$\nu = 2$}
\end{subfigure}
\begin{subfigure}{.32\textwidth}
\begin{center}
\includegraphics[width=1\textwidth]{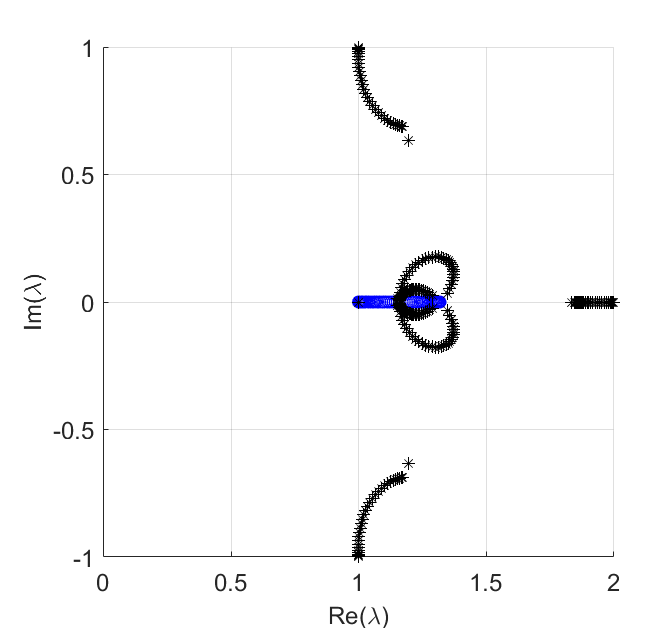}
\end{center}
\caption{$\nu = 4$}
\end{subfigure}
\begin{subfigure}{.32\textwidth}
\begin{center}
\includegraphics[width=1\textwidth]{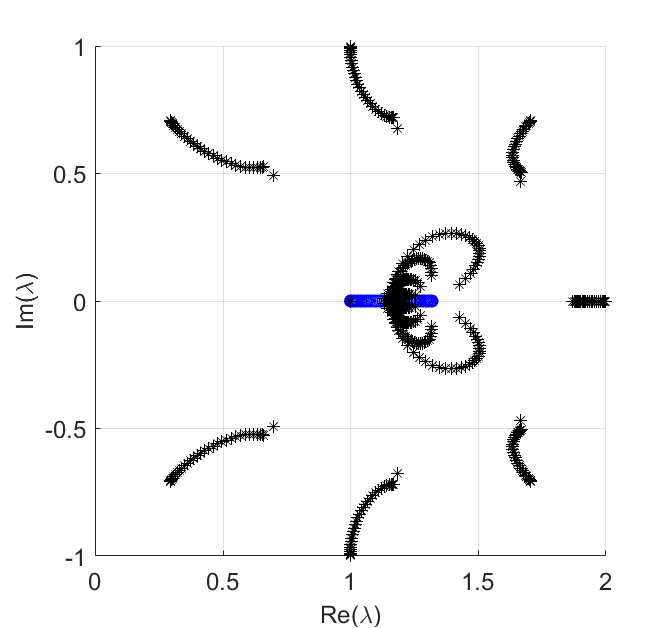}
\end{center}
\caption{$\nu = 8$}
\end{subfigure}
\end{center}
\caption{Spectrum of $\mathcal{A}_h \mathcal{Q}_h$ (in black), from 1-dimensional heat equation. The two-level setting uses time-space coarsening integrid operators~\eqref{eq:TS-coarsening}, ~\eqref{eq:constructZ} and~\eqref{eq:constructY}. The spectrum of $\mathcal{A}_h$ is plotted in blue. For the discretization, $n_x = 63$ and $n_t = 32$, corresponding to the Courant number $C = 2k\Delta t/\Delta x^2 = 0.32$, and $\theta = 0.5$.} \label{fig:spectrum_1d_aggloB}
\end{figure}

\begin{figure}[!h]
\begin{center}
\begin{subfigure}{.32\textwidth}
\begin{center}
\includegraphics[width=1\textwidth]{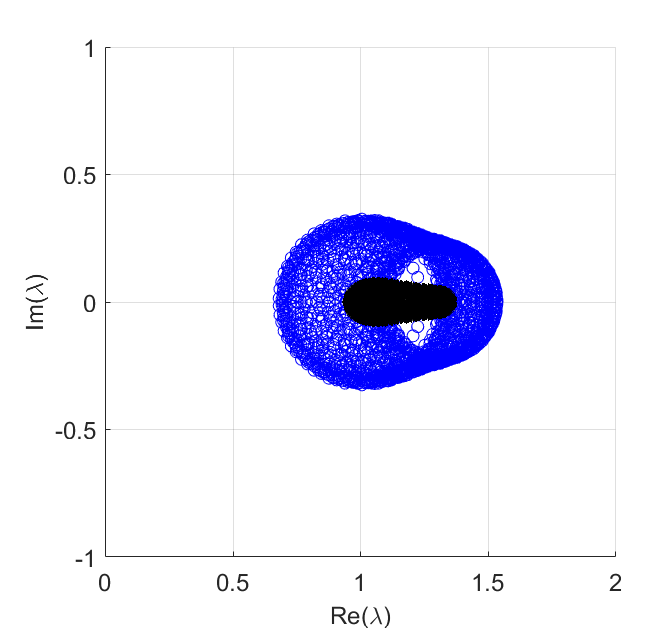}
\end{center}
\caption{$\nu = 2$}
\end{subfigure}
\begin{subfigure}{.32\textwidth}
\begin{center}
\includegraphics[width=1\textwidth]{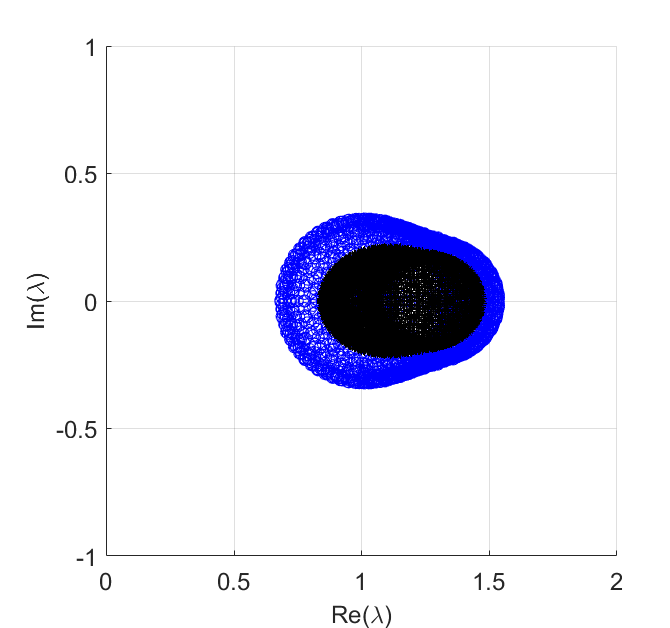}
\end{center}
\caption{$\nu = 4$}
\end{subfigure}
\begin{subfigure}{.32\textwidth}
\begin{center}
\includegraphics[width=1\textwidth]{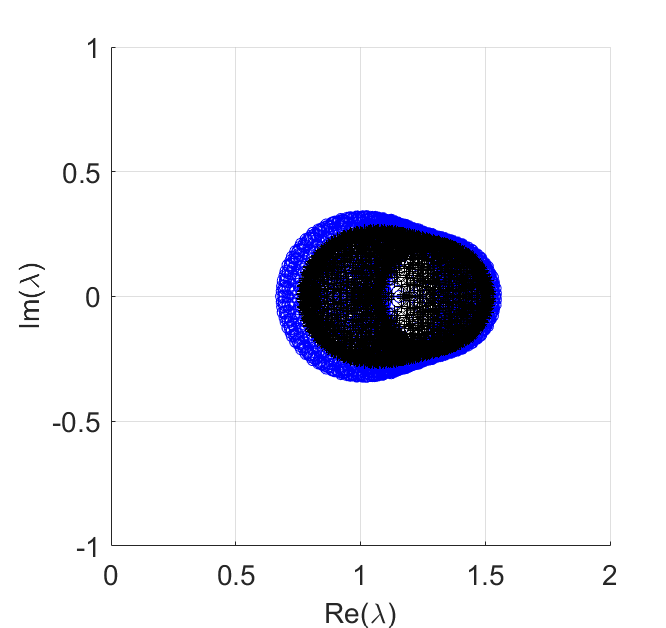}
\end{center}
\caption{$\nu = 8$}
\end{subfigure}
\end{center}
\caption{Spectrum of $\mathcal{A}_h \mathcal{Q}_h$ (in black), based on 1-dimensional heat equation. The two-level setting uses MGRIT intergrid operators~\eqref{eq:mgritY} and~\eqref{eq:mgritZ}, and the coarse-grid matrix~\eqref{eq:coarsemat}. The spectrum of $\mathcal{A}_h$ is plotted in blue. For the discrerization, $n_x = 64$ and $n_t = 32$, corresponding to the Courant number $C = 2\tau \Delta t/\Delta x^2 = 0.32$, and $\theta = 0.5$.} \label{fig:sprectrum_1d_mgrit}
\end{figure}

Spectra of the  two-level matrix with the time-space (TS) coarsening and MGRIT integrid components are shown in Figure~\ref{fig:spectrum_1d_aggloB}  and~\ref{fig:sprectrum_1d_mgrit}, respectively. For the TS-coarsening, the spectra are similar to the one with the T-coarsening, characterized by clusters of eigenvalues in the complex plane. Some eigenvalues are however spread near one. This cluster of eigenvalues tend to spread farther from one with an increase in $\nu$. The MGRIT-based two-level matrix leads to a completely different eigenvalue distribution from the T- and TS-coarsening. In this case, the eigenvalues lie inside a convex hull that contains the eigenvalues of the original matrix $\mathcal{A}_h$ (recall that in this case, the matrix $\mathcal{A}_h$ is rearranged according the fine-coarse partition).

With this type of spectra, it is not easy to deduce if a two-level Krylov approach will bring a significant convergence improvement, especially when $\nu > 2$, where the spectra become even much less clustered. Despite this spread of eigenvalues, the computed condition numbers $\kappa$ of the two-level matrix, presented in Table~\ref{tab:condition_number_1d}, are significantly lower than that of $\mathcal{A}_h$. Furthermore, the table suggests that a better conditioning can be achieved by the T-coarsening strategy. The T-coarsening also leads to the condition number of the two-level matrix which is practically independent of the Courant number $C = 2\tau \Delta t/\Delta x^2$. For all coarsening strategies, the conditioning deteriorates with increasing $\nu$, with the T-coarsening showing the least deterioration. This suggests that aggresive time coarsening (large $\nu$) may not be a good scenario for achieving a low number of iterations to reach convergence. For the TS-coarsening, because $\text{rank} Z < \text{rank} I_n = n$, we expect  less near-zero eigenvalues being shifted to 1, leading to slower convergence than with the T-coarsening.

\begin{table}[!h]
\caption{Condition numbers of $\mathcal{A}_h\mathcal{Q}_h$, based on the 1-dimensional heat equation. For the discretization, $n_x = 64$, with $\Delta t$ chosen according to the Courant number $C = 2 \tau \Delta t/\Delta x^2$, $\theta = 0.5$.} \label{tab:condition_number_1d}
\begin{center}
\begin{tabular}{|l|c|ccc|c|} \hline
                                 &  \multirow{2}{*}{$C = \frac{2\tau \Delta t}{\Delta x^2}$}    &   \multicolumn{3}{|c|}{$\kappa (\mathcal{A}_h \mathcal{Q}_h)$} & $\kappa(\mathcal{A}_h)$\\
                                 &  & $\nu = 2$ &  $\nu = 4$  & $\nu = 8$  &  \\ \hline
  \multirow{3}{*}{T-coarsening }  & 0.16 & 2.62 & 4.05 & 8.69 & 82.95 \\
                               & 0.32 & 2.62 & 4.05 & 8.69 & 42.39 \\
                               & 0.64 & 2.62 & 4.05 &  8.67 & 22.11 \\ \hline
  \multirow{3}{*}{TS-coarsening }  & 0.16 & 17.20 & 19.13 & 24.05 & 82.95 \\
                              & 0.32 & 8.82 & 9.85 & 12.87 & 42.39  \\
                              & 0.64 & 4.84  & 5.53  & 8.87&  22.11 \\ \hline
  \multirow{3}{*}{ MGRIT-based intergrid operators } & 0.16 & 2.64 & 4.05 & 10.76 & 82.95 \\
                             & 0.32 & 2.69 & 5.52 & 10.73 & 42.39  \\
                             & 0.64 & 2.85 & 5.59  & 10.74 &  22.11 \\ \hline
\end{tabular}
\end{center}
\end{table}

To see how the above results on condition numbers are translated to the actual convergence, we performed a numerical test by solving the two-level system~\eqref{eq:AQsystem} using GMRES, with the fine (level 1) time-space discretization corresponding to the Courant number $C = 2\tau \Delta t/\Delta x^2 = 0.64$. In this and all other numerical tests, we set $\bsu^0_h = 0$ and declare convergence at the $J$-th iteration if the following criteria is met:
$$
  \| \bsf_h - \mathcal{A}_h \bsu^J_h \|_2 / \| \bsf_h \|_2 < 10^{-6}.
$$
We note here in passing that in all numerical tests presented here, at convergence, the norm of difference between the solutions of the two- or multilevel Krylov method and the $\theta$-scheme is of order of $10^{-5}$.

The convergence results in terms of iteration counts are shown in Table~\ref{tab:1d_2lev_01}. As expected from its low condition numbers, the T-coarsening consistently leads to a more superior approach than the other two approaches. The plain use of MGRIT components does not generally lead to a better approach than the other options, except when the standard time coarsening ($\nu =2$) is used. For the latter case, MGRIT-based approach converges faster than the TS-coarsening approach.

\begin{table}[h!]
\caption{Numbers of GMRES iterations for solving the 1-dimensional heat equation with two-level Krylov: $C = 2\tau \Delta t/\Delta x^2 = 0.64$, $\theta = 0.5$} \label{tab:1d_2lev_01}
\begin{center}
\begin{tabular}{|c|c|ccc|ccc|ccc|} \hline
              &           &  \multicolumn{3}{c|}{T-coarsening} &  \multicolumn{3}{c|}{TS-coarsening}  &  \multicolumn{3}{c|}{MGRIT-based}\\
             &            &  \multicolumn{3}{c|}{$\nu:$} &  \multicolumn{3}{c|}{$\nu:$}  &  \multicolumn{3}{c|}{$\nu:$}\\ 
   $n_x$  & $n_t$  &    2    &   4   & 8 & 2  &   4   &   8  &  2  &  4  &  8 \\ \hline
   128     &   128   &  6  &  16  & 35 & 19 &  24 & 38   &   12    &  25  & 50\\
   256     &   512   &  6  &  15  & 32 & 18 &  22 & 36  &   11    &  23  & 49 \\
   512     &  2048  &  5  &  14  & 30 & 18 & 22 & 34 &   11    &  22  & 43 \\
  1024    &  8192  &  5  &  13  & 27 & 17 & 21  & 32  &   10    &  22  & 42 \\ \hline
\end{tabular}
\end{center}
\end{table}

%
%

\begin{figure}[!h]
\begin{center}
\includegraphics[width=0.6\textwidth]{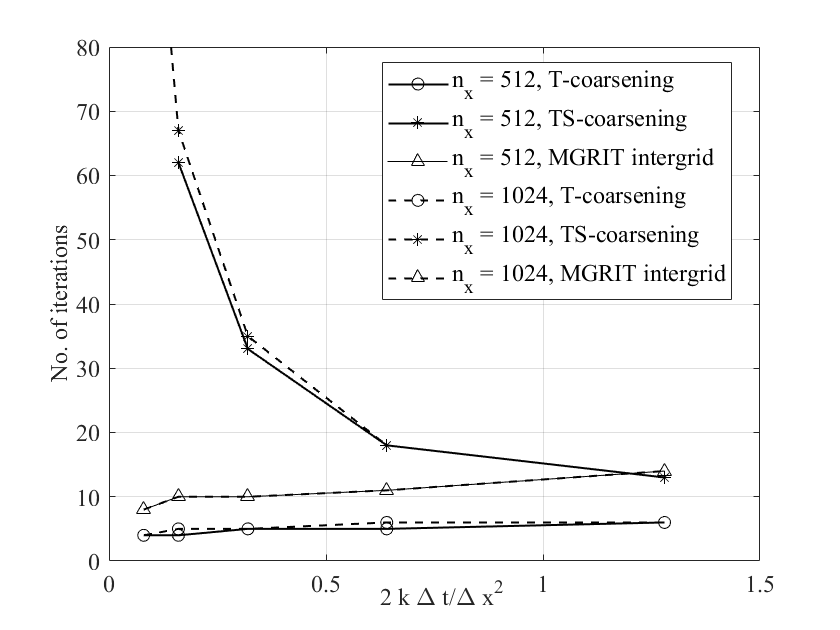}
\end{center}
\caption{The effect of the Courant number $C= 2\tau  \Delta t/\Delta x^2$ on the number of iterations.  The test problem is the 1-dimensional heat equation.} \label{fig:convergence_cfl_1d}
\end{figure}

Figure~\ref{fig:convergence_cfl_1d} shows iteration counts to reach convergence as a function of the Courant number $C$, with $n_x = 512$ and $1024$. Again, the T-coarsening approach results in the fastest two-level method, with iteration counts to convergence practically independent of the Courant number considered in this test. The MGRIT-based two-level method shows a slight increase in the iteration count as $C$ increases.  The TS-coarsening, while does not result in the fastest method in terms of iteration counts, shows a significantly improved convergence in cases where $C$ is high. As suggested by the increased condition numbers with larger $\nu$, aggressive coarsening affects the convergence negatively. 

\subsubsection{2-dimensional case} \label{subsec:2L-2d}

In this section, we discuss computational performance of the 2-level Krylov method for solving the 2-dimensional heat equation. We use the standard time coarsening $(\nu = 2)$ in the T-coarsening and MGRIT-based approach and additional standard space coarsening in the TS-coarsening. The computation is performed using Matlab on an Intel Core i7-7th Gen machine with 24 Gbytes of RAM. For the tests, the numbers of spatial gridpoints in the $x$- and $y$-direction, denoted by $n_x$ and $n_y$ respectively, are set to be equal, with $n_x \in \{128, 192\}$. Even though the method is implemented sequentially, all parallelizable parts in the method are coded in a way that mimics parallel computations. The  computational time thus is measured as if such a task is performed in parallel in multiple processors. This {\it simulated parallel} computation only gives us an approximation to the actual parallel computational time, as practically there is no communication time. 

\begin{figure}[!h]
\begin{center}
\begin{subfigure}{.49\textwidth}
\begin{center}
\includegraphics[width=1\textwidth]{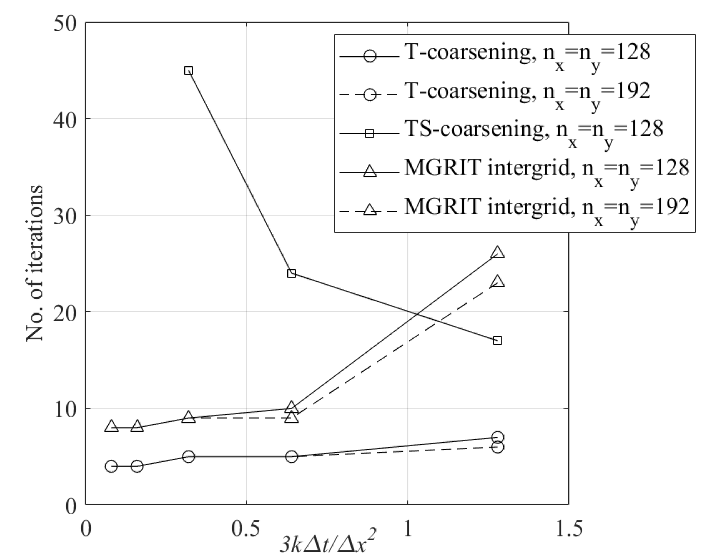}
\end{center}
\end{subfigure}
\begin{subfigure}{.395\textwidth}
\begin{center}
\includegraphics[width=1\textwidth]{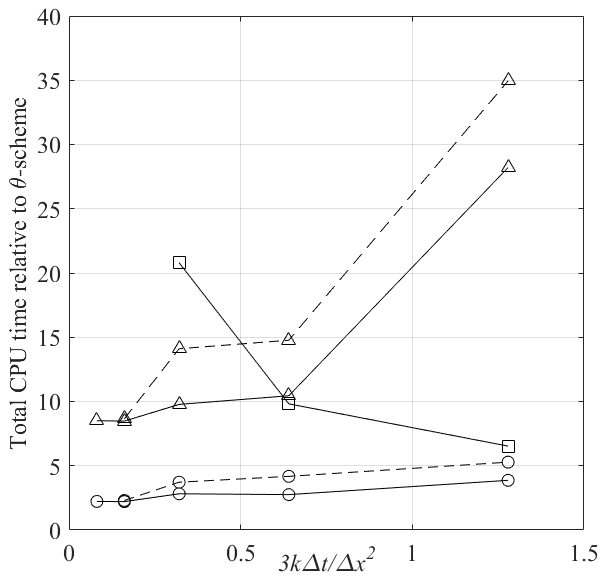}
\end{center}
\end{subfigure}
\begin{subfigure}{.405\textwidth}
\begin{center}
\includegraphics[width=1\textwidth]{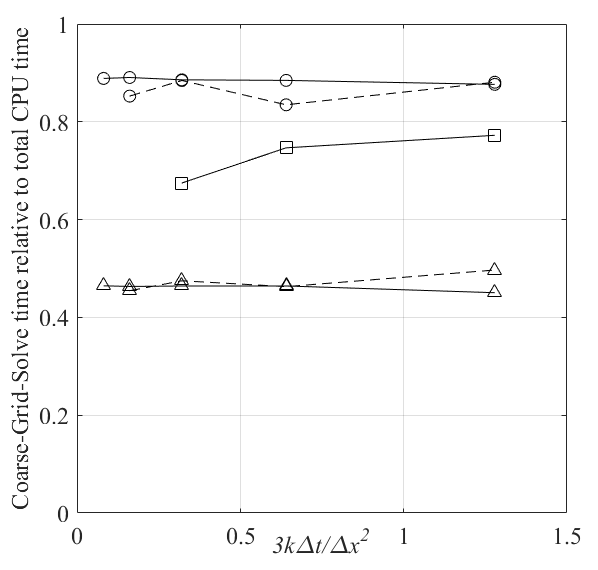}
\end{center}
\end{subfigure}
\end{center}
\caption{Convergence performance of the sequential two-level Krylov on single processor for solving the 2-dimensional heat equation with the T- and TS-coarsening, and MGRIT-based approach for various values $C = 3\tau \Delta t/\Delta x^2$. Top-left figure: numbers of iterations; Top-right figure: total time relative to the sequential $\theta$-scheme; Bottom figure: Coarse time-grid solve (CGS) time relative to total time.} \label{fig:CFL_Performance_2d}
\end{figure}

Convergence results of the sequential (single processor) two-level Krylov method in terms of numbers of iterations for varying values of $C =  3\tau \Delta t/\Delta x^2$ are shown in Figure~\ref{fig:CFL_Performance_2d}:top-left. As seen from the figure, the T-coarsening outperforms the other two approaches in terms of iteration counts. For the T-coarsening and MGRIT-based approach, the convergence rate deteriorates  as the Courant number $C$ increases, with the latter at a much higher rate. The convergence of the TS-coarsening approach improves when $C$  increases towards 1 but deteriorates as $C$ increases further (not shown in the figure). Overall, this convergence behavior is consistent with the behavior we observed earlier in the 1-dimensional case; cf. Section~\ref{sec:1d-case}.

The T-coarsening strategy also results in the fastest method in terms of the total computational time; see Figure~\ref{fig:CFL_Performance_2d}:top-right. Note that, in this figure, the time is presented as the time relative to the computational time of the sequential $\theta$-scheme. All methods, however, require more computational time than the $\theta$-scheme, with the T-coarsening two-level Krylov method be slower by a factor of at least 2. Due to the inversion needed in the interpolation step, the MGRIT-based two-level method becomes much slower than the T-coarsening. With $\nu = 2$, for instance, the MGRIT-based two-level method has to perform $n_t/2$ inversions of a matrix of size of $n^2$, with $n = n_x n_y$ in the interpolation step per iteration in addition to $n_t/2$ inversion of a matrix of the same size in the coarse time-grid solve. In constrast, the simple interpolation used in the T-coarsening does not involve inversion and costs only of order $n = n_x n_y$ computational work, which is negligible compared to the cost of inversion. Consequently, for the T-coarsening method, the computational work is dominated (about 90\%)  by the coarse time-grid solve (see Figure~\ref{fig:CFL_Performance_2d}:bottom).

\begin{figure}[!h]
\begin{center}
\includegraphics[width=0.45\textwidth]{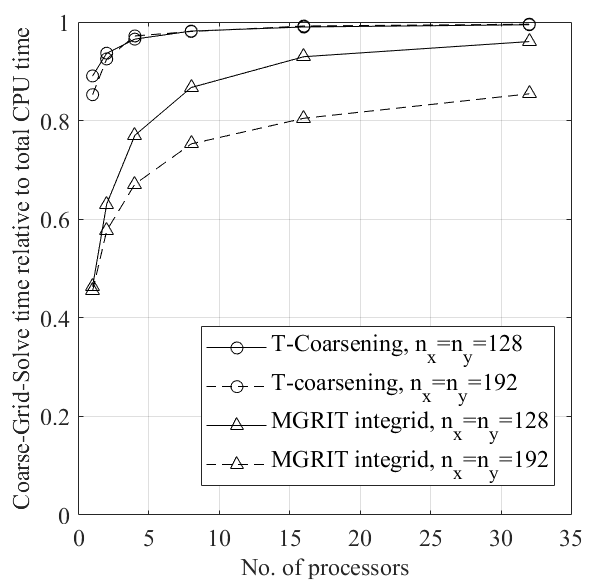}
\end{center}
\caption{Time spent on the coarse time-grid solve in the simulated parallel setting. The underlying problem is the 2-dimensional heat equation, discretized with $C = 3\tau\Delta t/\Delta x^2 = 0.16$.} \label{fig:twolevel_parallel_2D}
\end{figure}

Computational time can be reduced via parallel implementation of the method. Note that in the two-level approach, apart from the matrix inversion, the other operations such as matrix-vector multiplications, restriction, and interpolation can naturally be performed in parallel. Theoretically, taking into account only the matrix-vector multiplication and the coarse time-grid solve and assuming the inversion work of order of $n_x^2 n_y^2$ for the 2-dimensional heat equation, the work a GMRES-based two-level method requires to reach convergence is 
$$
  W_{\text{2-Level}} \sim  2 n_{\text{iter}}^{(1)} n_t n_x n_y + n_{\text{iter}}^{(1)} \frac{n_t}{\nu} n_x^2 n_y^2,
$$
with $n_{\text{iter}}^{(1)}$ the number of GMRES iteration on the first (finest) time level. Comparing with the $\theta$-scheme, whose total work is given by
$$
   W_{\text{$\theta$-scheme}} \sim n_t n_x n_y + n_t n_x^2 n_y^2,
$$
it is practically impossible to make a two-level multilevel Krylov faster than the $\theta$-scheme without parallelization. With parallelization of matrix-vector operations, the computational work on a processor  is given by
$$
   W_{\text{2-Level}} \sim  \frac{2 n_{\text{iter}}^{(1)}}{N_{\text{proc}}}  n_t n_x n_y + n_{\text{iter}}^{(1)} \frac{n_t}{\nu} n_x^2 n_y^2,
$$
where $N_{\text{Proc}}$ is the number of processors. In this case, the computational work on matrix-vector operations can be made much less than that of the $\theta$-scheme using a large number of processors. The effect of parallelization of the matrix/vector operations on the computational time is shown in Figure~\ref{fig:twolevel_parallel_2D}. As seen in the figure, parallel implementation of the T-coarsening approach reduces the computational time to a point that the contribution of the matrix/vector operations to the overall computational time becomes negligible. In the MGRIT-based approach, even though the interpolation can also be performed in parallel, matrix inversions involved in the interpolation process has a non-negligible contribution to the total computational time. 

Increased parallelism in the coarse time-grid solve can be attained by zeroing some off-diagonal parts of  the coarse time-grid matrix $\mathcal{A}_H$, hence decoupling the coarse time-grid system. For example, in the extreme case, we may drop all off-diagonal blocks and use the two-level Krylov matrix
\begin{eqnarray}
 \widetilde{\mathcal{Q}}_h = \mathcal{I}_h - \mathcal{Z} \widetilde{\mathcal{A}}_H^{-1} \mathcal{Y}^{\top} \mathcal{A}_h +\mu \mathcal{Z} \widetilde{\mathcal{A}}_H^{-1} \mathcal{Y}^{\top},   \label{eq:Qhtilde}
\end{eqnarray}
where $\widetilde{\mathcal{A}}_H = \text{blockdiag}(\mathcal{A}^H_{k,k})$. In this case, the coarse time-grid system is now completely decoupled and the solution procedure becomes parallelizable. For a choice of $N_\text{CGS}$ off-diagonal blocks set to zero, the course-grid system is decoupled into $N_\text{CGS}$ systems, which can be solved in parallel. In this case, with the approximate coarse-grid matrx $\widetilde{\mathcal{A}}_H$, one can show that 
\begin{align}
   \mathcal{A}_h\widetilde{\mathcal{Q}}_h = \begin{bmatrix}
                                                    \mathcal{A}^h_{0,0} \mathcal{Q}^h_0 &                &  & \\
                                                         \star   & \mathcal{A}^h_{1,1} \mathcal{Q}^h_1 &  & \\
                                                         \star   & \star         & \ddots  &  \\
                                                                  & \star         & \star  & \mathcal{A}^h_{n_T,n_T} \mathcal{Q}^h_{n_T} 
                                               \end{bmatrix}.
\end{align}
Hence, $\text{blockdiag}(\mathcal{A}_h\widetilde{\mathcal{Q}}_h) = \text{blockdiag}(\mathcal{A}_h \mathcal{Q}_h) $, and consequently,
\begin{align}
  \sigma(\mathcal{A}_h\widetilde{\mathcal{Q}}_h) = \sigma(\mathcal{A}_h \mathcal{Q}_h).
\end{align}

\begin{figure}[!h]
\begin{center}
\begin{subfigure}{.32\textwidth}
\begin{center}
\includegraphics[width=1\textwidth]{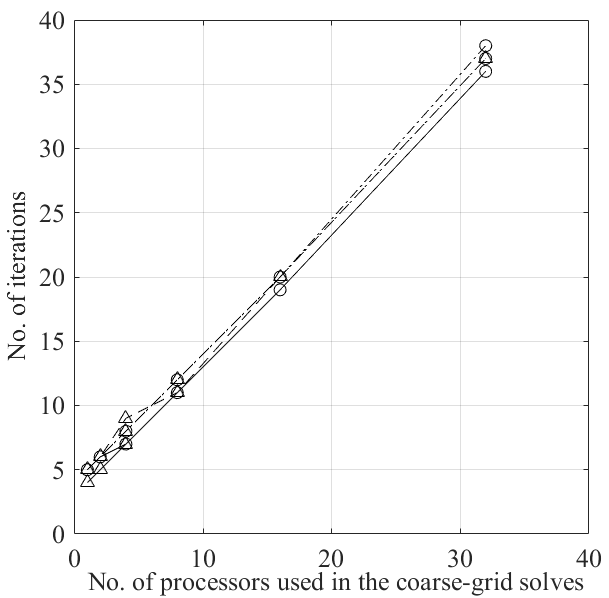}
\end{center}
\end{subfigure}
\begin{subfigure}{.32\textwidth}
\begin{center}
\includegraphics[width=1\textwidth]{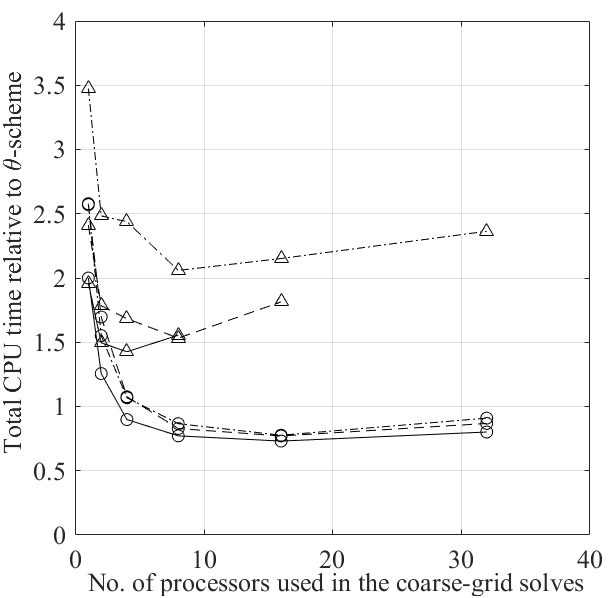}
\end{center}
\end{subfigure}
\begin{subfigure}{.325\textwidth}
\begin{center}
\includegraphics[width=1\textwidth]{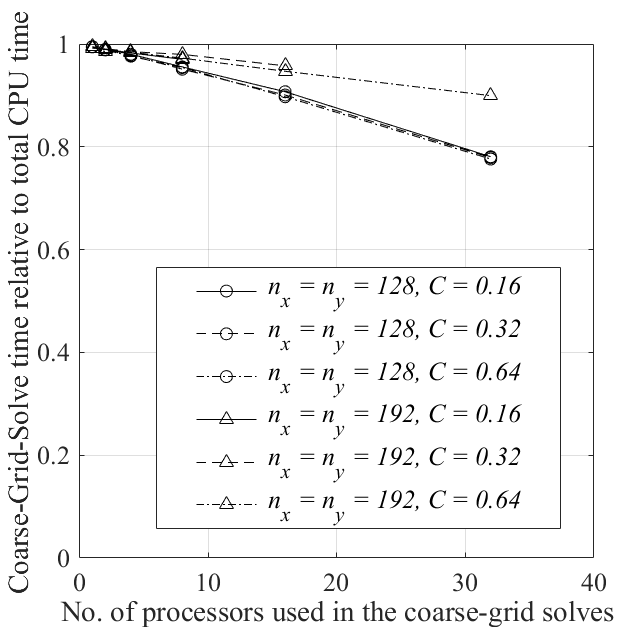}
\end{center}
\end{subfigure}
\end{center}
\caption{Effect of parallel coarse time-grid solves on the performance of two-level Krylov method with the T-coarsening and approximate coarse time-grid matrix $\widetilde{A}_H$. All matrix-vector operations are performed in parallel using 32 simulated  processors. The coarse time-grid system is decoupled into $N_{\text{CGS}} = N_{\text{proc}}$ systems. Left figure: iteration numbers; Middle figure: total time relative to the $\theta$-scheme; Right-figure: Time spent on the coarse time-grid solve relative to the total computational time.} \label{fig:twolevel_cgs_parallel_2D}
\end{figure}

Figure~\ref{fig:twolevel_cgs_parallel_2D}:left shows numbers of iterations to reach convergence using the T-coarsening and the approximate coarse time-grid matrix $\widetilde{\mathcal{A}}_H$. All matrix-vector operations are performed in parallel using 32 simulated processors. The coarse time-grid system is decoupled into $N_{\text{CGS}}$ systems by zeroing $N_{\text{CGS}} = N_{\text{Proc}}$ off-diagonal blocks, such that each simulated processor has an equal load, equivalent to inverting an $\frac{n_T}{N_{\text{CGS}}} \times \frac{n_T}{N_{\text{CGS}}}$ block matrix. As shown in Figure~\ref{fig:twolevel_cgs_parallel_2D}:left, the number of iterations increases linearly in the number of processors, but is almost indendent of $C$. Thus, even though $\mathcal{Q}_h$ and $\widetilde{Q}_h$ lead to spectrally equivalent two-level Krylov matrices, this fact does not necessarily translate to an idential convergence behavior, because eigenvalues are not the only factor that determines the convergence. The parallel coarse-grid solve however reduces the total computational time down significantly, even though there is a tendency to increase again with a larger number of used processors. In particular, with $n_x = n_y = 192$, the two-level Krylov method outperforms the $\theta$-scheme (Figure~\ref{fig:twolevel_cgs_parallel_2D}:mid). We note here that the total work with full parallelization of the two-level method is
$$
  W^{\text{full parallel}}_{\text{2-level}} \sim \frac{2 n_{\text{iter}}^{(1)}}{N_{\text{proc}}}  n_t n_x n_y +\frac{ n_{\text{iter}}^{(1)}}{\nu N_{\text{CGS}}} n_t n_x^2 n_y^2.
$$
With $N_{\text{CGS}} = N_{\text{proc}}$ and  $n^{(1)}_{\text{iter}}/N_{\text{CGS}} < 1.5$ as observed in Figure~\ref{fig:twolevel_cgs_parallel_2D}:left,  we have the constant $ n^{(1)}_{\text{iter}}/(\nu N_{\text{proc}}) < 1$, justifying this observed speed-up.  The propotion of work on the coarse time-grid solve also decreases (Figure~\ref{fig:twolevel_cgs_parallel_2D}:right; compare this with Figure~\ref{fig:twolevel_parallel_2D}) when more processors are used.


\section{Parallel-in-Time Multilevel Krylov and numerical results}  \label{sec:multilevel}

In this section, we discuss the extension of the two-level Krylov method discussed in Section~\ref{sec:two-level} to the multilevel Krylov method. We note that the multilevel Krylov method is flexible in the way the coarse-grid system is solved. What is crucial in the multilevel Krylov method is the accuracy and the effeciency of the coarse-grid solve. One can employ, e.g., MGRIT, to solve the coarse-grid system. A standard multilevel Krylov method, however, is built upon a recursive application of the two-level Krylov method. For elliptic PDEs, as the spatial grid becomes coarser, the underlying coarse space-grid matrix tends to become better conditioned. Hence, the number of iterations needed to approximately solve the coarse-grid system on the coarser levels can be reduced. As reported in~\cite{Erlangga08Na, Erlangga08Nc}, a multilevel Krylov method needs only a few GMRES iterations on the coarser levels. 

\subsection{Time coarsening}

For the T-coarsening multilevel Krylov method with $L$ levels, it can be shown the computational work is
\begin{align}
  W^{\text{full parallel}}_{\text{$L$-level}} \sim \frac{2n_{\text{iter}}^{(1)}}{N_{\text{proc}}} \left( 1 + \sum_{\ell = 2}^{L-1} \frac{1}{\nu^{\ell-1}} \prod_{j=2}^{\ell} n_{\text{iter}}^{(j)} \right)n_t n_x n_y + \frac{n_{\text{iter}}^{(1)}}{N_{\text{CGS}}} \left(\prod_{\ell = 2}^{L-1} \frac{n_{\text{iter}}^{(\ell)}}{\nu^{L-2}} \right) n_t n_x^2 n_y^2. \label{eq:West_MK}
\end{align}
The above estimate suggests that the multilevel implementation adds extra cost in matrix/vector operations. The cost of coarse time-grid solve can however potentially be reduced if the number of iterations used on the coarse levels is less than $\nu$. For standard coarsening, this will require $n_{\text{iter}}^{(\ell)} \le \nu = 2$. As observed in the experiments, this requirement is quite stringent if applied on every level, and results in slow convergence due to inaccuracy in the coarse-grid solve. Because the accuracy of the second-level coarse-grid solve seems to be crucial, we shall allow more GMRES iterations on the second level and compensate the extra cost they bring by performing only one or two iterations on the coarser levels.

\begin{figure}[!h]
\begin{center}
\begin{subfigure}{.32\textwidth}
\begin{center}
\includegraphics[width=1\textwidth]{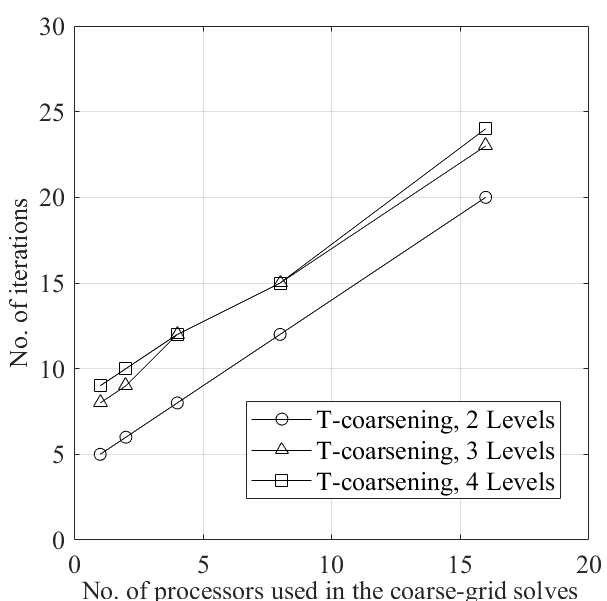}
\end{center}
\end{subfigure}
\begin{subfigure}{.32\textwidth}
\begin{center}
\includegraphics[width=1\textwidth]{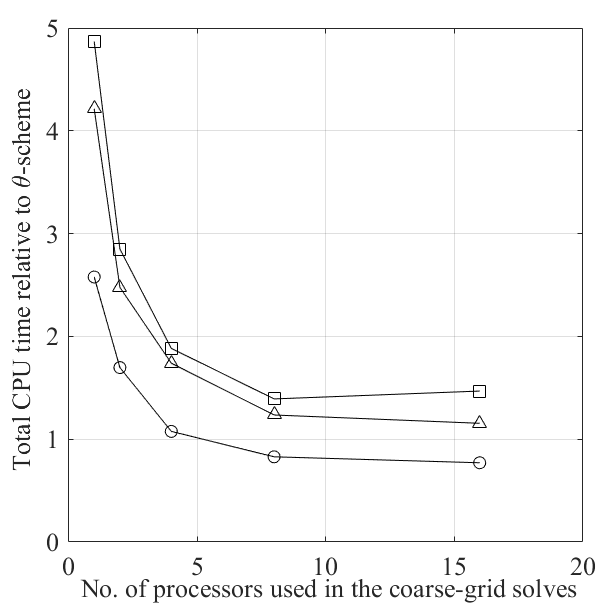}
\end{center}
\end{subfigure}
\begin{subfigure}{.325\textwidth}
\begin{center}
\includegraphics[width=1\textwidth]{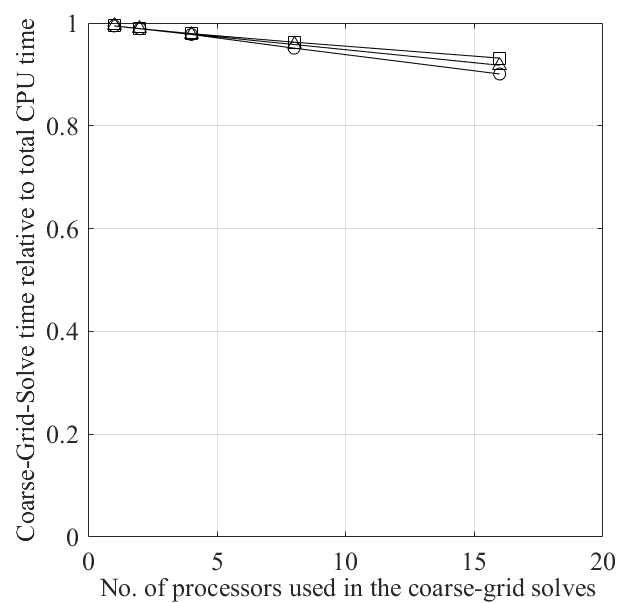}
\end{center}
\end{subfigure}
\end{center}
\caption{Performance comparison of multilevel Krylov with time coarsening and various numbers of processors used to solve the coarse-grid system in parallel, with $n_x = n_y = 128$, $n_t = 768$. Other than the coarse grid solve, 32 processors are used to perform all related matrix/vector operations.} \label{fig:MK-t-coarsening-032}
\end{figure}

Figure~\ref{fig:MK-t-coarsening-032} shows performance of 2-, 3-, 4-level Krylov method with T-coarsening and with standard coarsening ($\nu = 2$) applied to generate a sequence of coarse time-grid systems. Matrix-vector operations at all levels are performed on 32 simulated processors. On the coarsest level, the coarse time-grid system is solved in parallel using the perturbed coarse-grid matrix $\widetilde{A}_H$, with $N_{\text{CGS}} = N_{\text{proc}}$. At any level between the finest and the coarsest level, 2 GMRES iterations are used. 

As shown in Figure~\ref{fig:MK-t-coarsening-032}:left, with T-coarsening, the number of iterations tends to increase as more  coarse levels are used. This in turn increases the total CPU time needed to reach convergence; see Figure~\ref{fig:MK-t-coarsening-032}:mid. In this case, the parallel multilevel approach does not lead to a faster method than the sequential $\theta$-scheme. Even though more iteration numbers have to be spent when more processors are used on the coarsest level, the smaller coarse time-grid system each processor has to solve results in reduced CPU times. As shown in Figure~\ref{fig:MK-t-coarsening-032}:right, the coarse time-grid solve remains the most time consuming operation, taking more than 90\% of the total time.

We remark here that for time-dependent problems, coarser time grid may cause instability, leading to highly inaccurate solve of the coarse time-grid system. This inaccurate coarse time-grid solve can probably be compensated by allowing more GMRES iterations on the fine levels. This however may in the end render the method inefficient.

\subsection{Alternate time-space coarsening}

Based on the estimate~\eqref{eq:West_MK}, the computational time can be reduced further via the inclusion of spatial coarsening. For the 2d case, the use of standard coarsening in space (roughly doubling the spatial mesh size) at a level reduces the first and the second term in the estimate by a factor of 2 and 4, respectively.  One strategy for spatial coarsening is by performing time-space coarsening simultaneously. However, as suggested by the test results in Section~\ref{sec:1d-case}, this approach leads to a slow-to-converge two-level Krylov method. As the coarse-grid solve needs to be done as accurately as possible, applying a time-space coarsening at a particular level will require a significant number of iterations on the level the simultaneous time-space coarsening is performed. 

\begin{figure}[!h]
\begin{center}
\includegraphics[width=0.45\textwidth]{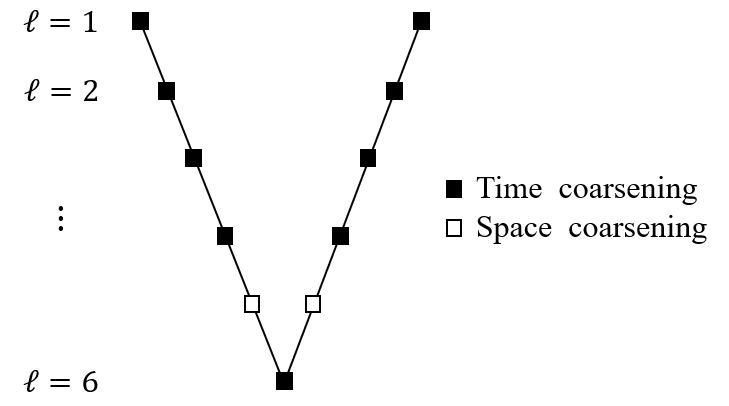}
\end{center}
\caption{An example of time-space coarsening, resulting in 6 levels.} \label{fig:MK-TS-coarsening}
\end{figure}

An alternative to the simultaneous time-space coarsening is by applying spatial coarsening on some coarse levels whenever the time grid has been sufficiently coarsened. This strategy is illustrated in Figure~\ref{fig:MK-TS-coarsening}, which involves 5 coarsening. Starting from the fine time-space grid ($\ell = 1$), the next three coarse-level systems are constructed via time coarsening, indicated by the black squares. After this sequence of time coarsening, the remaining coarse-level systems are constructed via space and time coarsening, applied alternatingly (the space coarsening is indicated by the white square). This coarsening strategy results in a coarse-level system associated with $n_t/8$ time steps $(\Delta T = 8\Delta t)$. In each time step, a coarse matrix of size $n_x n_y/4 \times n_x n_y/4$  has to be inverted.  For the numerical results presented here, the first space coarsening is performed whenever the sequence of time coarsening in the finer level results in time-space grids which violates the stability condition.

\begin{figure}[!h]
\begin{center}
\begin{subfigure}{.32\textwidth}
\begin{center}
\includegraphics[width=1\textwidth]{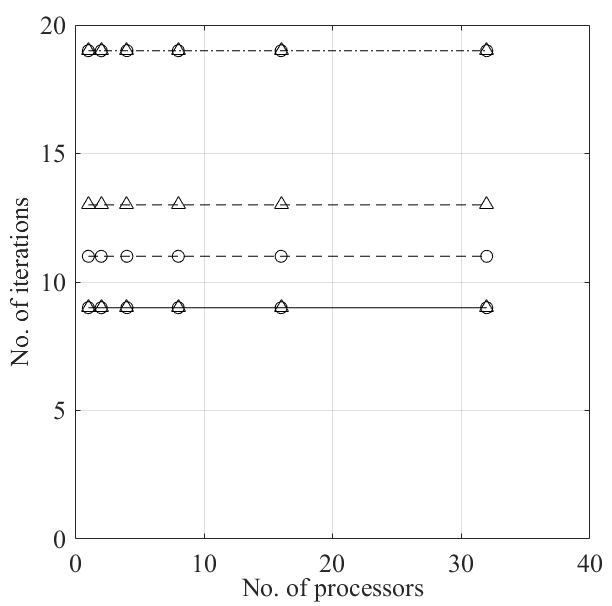}
\end{center}
\end{subfigure}
\begin{subfigure}{.32\textwidth}
\begin{center}
\includegraphics[width=1\textwidth]{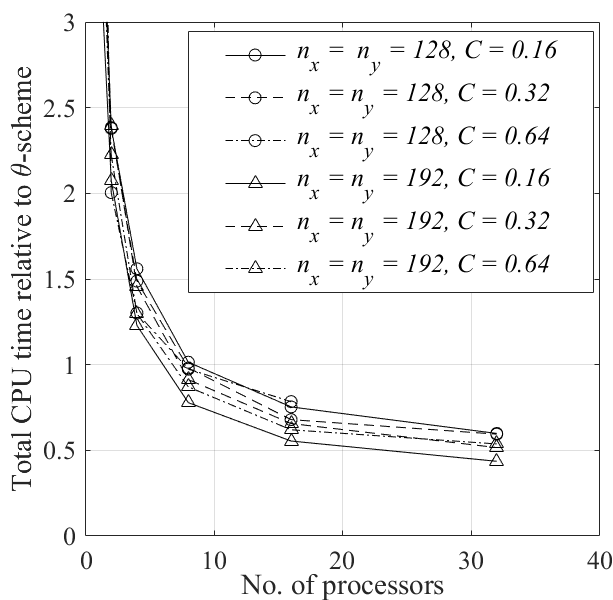}
\end{center}
\end{subfigure}
\begin{subfigure}{.325\textwidth}
\begin{center}
\includegraphics[width=1\textwidth]{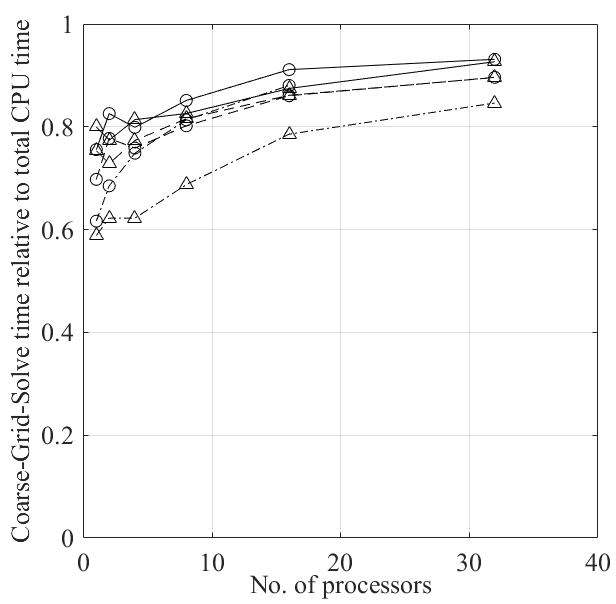}
\end{center}
\end{subfigure}
\end{center}
\caption{Performance comparison of multilevel Krylov with time-space coarsening and various numbers of processors used to solve the coarse-grid system in parallel. On the coarsest level, 16 processors is used to solve the coarse-grid system.} \label{fig:MK-ts-coarsening} 
\end{figure}

We show computational performance of the multilevel Krylov method with the time-space coarsening outlined above in Figure~\ref{fig:MK-ts-coarsening}. The computations were performed for various Courant numbers $C$, with two spatial gridpoint settings: $n_x = n_y = 128$ and $192$. In this case, with $C = 0.16$ for example, the initial time coarsening is performed until the 4th coarse level, followed by alternating space-time coarsening. On each level, computations are done in parallel using a fixed number of simulated processors, except on the coarsest level, where in all cases, the coarsest-grid system is decoupled into 16 systems, solved in parallel on 16 simulated processors. On the coarse levels, 2 FGMRES iterations are performed, except on the level before the coarsest one, where 1 FGMRES iteration is used.

As shown in Figure~\ref{fig:MK-ts-coarsening}:left, iteration numbers to converge depend on the Courant number $C$, with lower $C$ resulting in a faster convergence. The results also suggests that the number of iterations to converge is mildly dependent on the grid settings, so long as $C$ is kept constant. The total CPU time relative to the sequential $\theta$-scheme is shown in Figure~\ref{fig:MK-ts-coarsening}:mid. The graph demonstrates the gain in CPU time when a large number of processors is used (in this case, $n_{\text{proc}} > 8$), outperforming the serial $\theta$-scheme for all $C$ used. Furthemore, the computational time on the coarse-grid solve is reduced to below 90\% of the total computational time.

Finally, to investigate the effect of number of levels used in the method, we consider a case with $n_t = 2^{11}$ and $n_x = n_y = 192$ on the finest level, associated with $C = 0.16$. The time-space grid is coarsened in the similar way as in the previous test, with the coarsest level system solved in a decoupled way on 16 simulated processors. Apart from the coarsest-grid solves, matrix/vector operations are performed in parallel using either $32$ or $64$ simulated processors. The convergence results are shown in Table~\ref{tab:varlevel}. We note here that for the case with  32 processors and 7 levels, the data in the table corresponds to the right-end of the curve in Figure~\ref{fig:MK-ts-coarsening}, with $n_x = n_y = 192$ and $C = 0.16$. Thus, Table~\ref{tab:varlevel} suggests that the CPU time relative to the $\theta$-scheme can be further reduced via the use of more levels to decrease the size of the coarsest-grid system, without practically affecting the number of iterations. In this case, about 90\% of the computational time is spent on the coarse-grid solves. With 11 levels, the method achieves a computational speed-up of about 9.39 from the single-processor implementation. 

Using 64 simulated processors on the matrix/vector operations reduces the total CPU time further, as in this case, the computational time on the coarsest-grid solves increases to around 95\%; see Table~\ref{tab:varlevel}. Under this setting, we however stop the coarsening when reaching 64 time grids and after performing the necessary spatial coarsening (on Level 9). It is indeed possible that after Level 9, the coarsening is continued and the resulting coarser-grid systems are handled on, for instance, 32 simulated processors. This implementation scenario is, however, not considered here. Even with 9 levels, the method gains a computational speed up of about 9.85 relative to the single-processor implementation.  

\begin{table}[!h]
\caption{Performance comparison using different number of levels, using an alternate TS-coarsening, $n_x = n_y = 192$, $n_t = 2^{11}$. On the the coarsest level, 16 processors are used. Other than the coarsest level, all matrix/vectors operations are performed on 32 or 64 processors. The symbol ``--'' means coarser levels that can be fit in the assigned number of processors could not be constructed further.} \label{tab:varlevel}
\begin{center}
\begin{tabular}{|c|c|ccc|ccc|} \hline
          & coarsest level & \multicolumn{6}{|c|}{No. of Processors} \\ \cline{3-8}
          & time-space & \multicolumn{3}{|c|}{32} & \multicolumn{3}{|c|}{64} \\  \cline{3-8}
  Level & grid  & \#Iter. & $T_{\text{CPU}}$  &  $T_{\text{CGS}}$ &\#Iter. & $T_{\text{CPU}}$  & $T_{\text{CGS}}$ \\ \hline
    7    & $128 \times 48 \times 48$  & 9    & 0.436 & 0.927 & 10 & 0.322 & 0.948 \\
    8    & $64 \times 48 \times 48$  & 10   & 0.428 & 0.913 & 10 & 0.365 & 0.947 \\
    9    & $64 \times 24 \times 24$ &10   & 0.416 & 0.918 & 10 & 0.365 & 0.951 \\
   10   & $32 \times 24 \times 24$ & 10   & 0.422 & 0.916 &  --  &   --     &  --      \\
   11   &  $32 \times 12 \times 12$ & 9    & 0.384 & 0.916 &  --  &   --     &  --      \\ \hline
\end{tabular}
\end{center}
\end{table}

\section{Conclusion}
\label{sec:conclusion}

This paper discussed and extended the multilevel Krylov method for stationary problems to the time-dependent problems, which allows parallel-in-time computational strategies. Likewise parareal- and MGRIT-based methods, the method is iterative in nature and involves a coarse-level problem on which parallelization can be exploited in a multilevel fashion. As the coarsest level is typically the bottleneck for full time parallelization, we proposed an approximation to the coarsest-grid system by dropping some off-diagonal blocks in the coarsest-grid matrix, effectively decoupling the systems and allowing full time parallelization on all levels. This approximation was observed to be quite effective, especially when a sufficient number of coarse time levels was used.

Numerical tests suggested that a significant gain in computational time could not be achieved solely with time coarsening. The proposed method incorporated the space coarsening in an alternating way, following a sequence of first few time coarsening until a stability condition (dictated by the Courant number) is violated. This approached led to the total computational time as low as 0.32 relative to the sequential $\theta$-scheme. 

The current implementation was limited by the author's computational resources. It is therefore expected to perform an actual test of the method in the future on a parallel computer. Furthermore, as this paper was meant to present only the prototypical implementation of the method, future's work will include, among others, the effective construction of the deflation (prolongation and restriction) matrices, even though the current simple construction seemed to work very well, time-space coarsening and as well the multilevel strategies, which would lead to further gain in computational time.

\bibliographystyle{siamplain}
\bibliography{mybiblio}

\end{document}